%% file: article.tex
\documentclass[12pt,a4paper]{article}
\usepackage{import}

\import{}{declarations.tex}

\usepackage{latexsym} 

\usepackage{color}

\begin{document}

\title{The conformal measures of a normal subgroup of a cocompact Fuchsian group}
\author{Ofer Shwartz\\
\\
Weizmann Institute of Science,\ Rehovot, Israel\\
\\\\}
\maketitle

%+Abstract
\begin{abstract}

In this paper, we study  the conformal measures of a normal subgroup of a cocompact Fuchsian group. In particular, we relate  the extremal conformal measures to the eigenmeasures of a suitable Ruelle operator. Using Ancona's theorem, adapted to the Ruelle operator setting, we show that if the group of deck transformations $G$ is hyperbolic then the extremal conformal measures and the hyperbolic boundary of $G$ coincide. We then interpret these results in terms of the asymptotic behavior of cutting sequences of geodesics on a regular cover of a compact hyperbolic surface. 
\end{abstract}

\import{}{introduction.tex}

\import{}{hyperbolic.tex}
%\import{}{zd_covers.tex}

\appendix
\import{}{ancona.tex}

\import{}{appendix.tex}

\bibliography{bibfile}
\bibliographystyle{plain}
%-Bibliography

\end{document}

%% file: declarations.tex
\usepackage{amssymb}
\usepackage{latexsym}
\usepackage{graphicx}
\usepackage{graphics}
\usepackage[normalem]{ulem}
\usepackage{amsmath}
\usepackage{amsthm}
\usepackage{amssymb}
\usepackage{amsfonts}
\usepackage{color,soul}
\usepackage{enumerate}
\usepackage{enumitem}
\usepackage{mathrsfs}

%\usepackage{hyperref}

%\usepackage[left=0.6in,right=0.6in, top=0.6in, bottom=1.6in]{geometry}
%\usepackage{setspace}
%\onehalfspacing

\usepackage{tikz}
\usetikzlibrary{matrix}
\usetikzlibrary{arrows}
\usetikzlibrary{shapes}

\newtheorem{thm}{Theorem}[section]
\numberwithin{thm}{section}

\newtheorem{prop}[thm]{Proposition}

\newtheorem{lemma}[thm]{Lemma}

\newtheorem{cor}[thm]{Corollary}
\newtheorem{remark}[thm]{Remark}

\newtheorem{defn}[thm]{Definition}

\newtheoremstyle{named}{}{}{\itshape}{}{\bfseries}{.}{.5em}{#1 \thmnote{#3}}
\theoremstyle{named}
\newtheorem*{numberedtheorem}{Theorem}

\def \conf {{\mathrm{Conf}}}

\DeclareMathOperator\confG{Conf(\Gamma, \delta)}
\DeclareMathOperator\extconf{ext(\confG)}

\def \T {{\bf T}}
\def \G {{\bf G}}

\def \I {{\bf I}}

%% file: introduction.tex
\section{Introduction}
Let  $\mathbb{D} = \{z\in \mathbb{C} : |z|<1\}$ be the open hyperbolic unit disc and let $\partial \mathbb{D}=\{z\in \mathbb{C} :\ |z|=1\}$.  Let $\Gamma$ be  a Fuchsian group (a discrete subgroup of M\"obius transformations) which preserves $\mathbb{D}$. We denote by $\delta(\Gamma)$ the critical exponent of $\Gamma$ (see definition in Section \ref{sec:HyperbolicCovers}). A finite measure $\mu$ on $\partial \mathbb{D}$ is said to be\textit{ $(\Gamma, \delta)$-conformal} if for every $\gamma\in \Gamma$,
\begin{equation*}
\frac{d(\mu\circ \gamma)}{d\mu} = |\gamma'|^\delta. 
\end{equation*}
We denote by $\confG$ the collection of $(\Gamma, \delta)$-conformal measures and by $\extconf$ the extremal points of $\confG$.

Conformal measures have many applications in hyperbolic geometry:
\begin{itemize}
\item 
\uline{\textit{Geodesic-flow-invariant measures:}}
If $\mu_1, \mu_2$ are two $(\Gamma, \delta)$-conformal measures then the measure $m(\xi^-, \xi^+, t)=\frac{d\mu_1(\xi^-)d\mu_2(\xi^+)dt}{||\xi^-- \xi^+||^{2\delta}}$  projects to a geodesic-flow-invariant measure on $T^1(\mathbb{D} / \Gamma)$ (the unit tangent bundle of $\mathbb{D} / \Gamma$), see \cite{bedford_1991}. However, not every geodesic-flow-invariant measure is of this form, for example measures which are supported on periodic orbits.
\item
\uline{\textit{Horocycle-flow-invariant measures:}} If $\mu$ is a $(\Gamma, \delta)$-conformal measure then the measure $dm (\xi, s,t) = e^{\delta s}d\mu(\xi)dsdt$ projects to a horocycle-flow-invariant measure on $T^1(\mathbb{D} / \Gamma)$. Moreover, if the underlying surface $\mathbb{D} / \Gamma$ is a \textit{tame surface} then every ergodic horocycle-flow-invariant Radon measure which is not supported on a single  horocycle is of this form, see \cite{sarig_2010_tame}. Recently, Landesberg and Lindenstrauss derived a similar decomposition for Radon horospherical-flow-invariant measures in higher dimensions \cite{landesberg_2019}.
\item
\uline{\textit{Eigenfunctions of the Laplacian-Beltrami operator:}} If $\mu$ is a $(\Gamma, \delta)$-conformal measure and $P:\mathbb{D} \times \partial \mathbb{D}\rightarrow \mathbb{R}^+$ is the Poisson kernel, $P(z,\xi)=\frac{1-|z|^2}{|\xi - z|^2}$, then
$$h(z) = \int_{\partial \mathbb{D}}P(z,\xi)^{\delta}d\mu(\xi), \qquad z\in \mathbb{D}$$
is a positive $\Gamma$-invariant $\delta(\delta-1)$-eigenfunction of the Laplacian-Beltrami operator and every such eigenfunction arises in that way, see \cite{babillot_2004,karpelevic_1965}.
\end{itemize}

\paragraph{Known results on the classification of the conformal measures.}
The existence of a $(\Gamma, \delta)$-conformal measure was first proven
by Patterson \cite{patterson_1976} for the critical value $\delta = \delta(\Gamma)$
and by Sullivan in higher dimensions \cite{sullivan_1979}.  Later on, in
\cite{sullivan_1987} Sullivan showed that for non-cocompact groups with no
parabolic elements, a $(\Gamma, \delta)$-conformal measure exists iff $\delta
\geq \delta(\Gamma)$. In \cite{roblin_2011}, Roblin studied the conformal
measures in more general settings via a Martin boundary approach. A more
general class of measures,   \textit{quasiconformal measures}, has  been
considered  as well, see \cite{coornaert_1993,blachere_2011}.

 Furstenberg
\cite{furstenberg_1973} showed that if $\Gamma$ is \textit{cocompact}, namely $\mathbb{D} / \Gamma$ is compact, then the  Lebesgue measure is the unique  $(\Gamma, \delta(\Gamma))$-conformal measure and  there are no other  $(\Gamma, \delta)$-conformal measures for all $\delta>\delta(\Gamma)$. Variants of this result were   proven  by Dani \cite{dani_1978} for cofinite groups  and by  Burger \cite{burger_1990}  for geometrically finite groups. Their original motivation was the classification of the horocycle-flow-invariant measures.

If  $\mathbb{D} / \Gamma$ is a regular cover  of a compact hyperbolic surface and with nilpotent group deck transformations $G$, then there is a bijection between the set of all extremal $(\Gamma, \delta)$-conformal measures (for all $\delta \geq \delta(\Gamma)$) and the set of all homomorphisms from $G$ to $\mathbb{R}$, see \cite{lyons_1984} and also \cite{ledrappier_2007}.

In \cite{sarig_2008}, Schapira and Sarig studied the  horocycle-flow-invariant measures on $\mathbb{Z}^d$-covers (namely $G = \mathbb{Z}^d$) in terms of the almost-surely asymptotic direction of geodesics.

In \cite{kaimanovich_2000_horocycle}, Kaimanovich characterized
the ergodicity of the horocycle flow with respect to the Liouville measure,
namely the specific case where $\delta =1$.  See also \cite{roblin_2003} for an extensive study of the conformal measures in  negatively curved geometrically finite
manifolds. 

\paragraph{Conformal measures and eigenmeasures of the Ruelle operator.}
In this work we study the conformal measures of a normal subgroup of a cocompact Fuchsian group, namely  under the assumption that there exists a cocompact Fuchsian group $\Gamma_0$ with $\Gamma \vartriangleleft  \Gamma_{0}$.

 For such $\Gamma$, we show that   for every $\delta \geq \Gamma(\delta)$ there is  a linear $1-1$ correspondence between the extremal $(\Gamma, \delta)$-conformal measures and eigenmeasures of a suitable Ruelle operator. 

This correspondence is stated using the Bowen-Series coding. In more details, for a  co-compact Fuchsian group $\Gamma_0$ let $F_0\subseteq \mathbb{D}$ be a  fundamental domain for $\mathbb{D} / \Gamma_0$. In  \cite{bowen_1979}, Bowen and Series constructed (w.r.t. $F_0)$  a finite partition $\{I_a\}_{a\in S_0}$ of $\partial \mathbb{D}$ into closed arcs with disjoint interiors and a  finite set $\{e_a\}_{a\in S_0}\subseteq\Gamma_0$ s.t. the set $\{e_a\}_{a\in S_0}$ generates $\Gamma_0$ and the \textit{Bowen-Series map} $f_{\Gamma_0}:\partial \mathbb{D}\rightarrow \partial\mathbb{D}$, 
$$f_{\Gamma_0}(\xi ) = e_{a}^{-1}\xi, \quad \xi \in int( I_a) $$
induces a Markov partition of $\partial \mathbb{D}$, namely the space
$$\Sigma: = \{(\sigma_i) : \forall i\geq 0, \;\sigma_i \in S_0 \text{ and }int(f_{\Gamma_0}(I_{\sigma_i}))\cap int(I_{\sigma_{i+1}})\neq \varnothing \}$$  
along with the left-shift transformation is a subshift of finite type.
 Let $\pi_{\Sigma}:\Sigma\rightarrow\partial\mathbb{D}$  be the canonical projection,  $\pi_\Sigma(\sigma) \in \cap_{n\geq 0}{f_{\Gamma_0}^{-n}I_{\sigma_n}}$ (the intersection is a singleton, see \cite{bedford_1991}).
 For several other
important properties of the Bowen-Series coding, see Section \ref{sec:bowen_series_coding}.

Let $(X,T)$ be the \textit{group extension }of $\Sigma$ with $G=\Gamma_0 / \Gamma$, 
$$X = \biggl\{\bigl((\sigma_0, \gamma_0\Gamma),(\sigma_1, \gamma_1\Gamma) ,\dots\bigr) : (\sigma_i)\in \Sigma; \;\;\forall i\geq 0, \;\gamma_i \Gamma\in G \text{ and }\gamma_{i+1}\Gamma = e_{\sigma_i}^{-1}\gamma_i \Gamma\biggr\} $$
and let $T:X\rightarrow X$ be the left-shift transformation, see  \cite{stadlbauer_2013}.
We sometime  use the following canonical correspondence to identify $X$ with $\Sigma \times G$,
$$(\sigma, \gamma\Gamma) \longmapsto \bigl((\sigma_0, \gamma\Gamma),(\sigma_1, e_{\sigma_0}^{-1}\gamma\Gamma),(\sigma_2, e_{\sigma_1}^{-1}e_{\sigma_0}^{-1}\gamma\Gamma),\dots\bigr).$$ 
Given $\delta>0$, let $\phi^{
X,\delta}:X\rightarrow\mathbb{R}$,  
$$\phi^{X,\delta}(\sigma ,\gamma\Gamma) := -\delta_{}\log|(e_{\sigma_0}^{-1})'(\pi_{\Sigma}(\sigma))| .$$
The \textit{Ruelle operator}, evaluated on a function $f:X\rightarrow\mathbb{R}$ and a point $x\in X$ is
$$(L_{\phi^{X,\delta}}f)(x) =  \sum_{y:Ty = x}e^{\phi^{X,\delta}(y)}f(y). $$
See Definition \ref{def:ruelle_operator}.
In Section \ref{sec:conf_bijection}, we prove the following theorem which connects between the conformal measures and the eigenmeasures of $L_{\phi^{X,\delta}}$.

\begin{numberedtheorem}[\ref{prop:gamma_conf_bijection}]
Let $\Gamma_0$ be a cocompact Fuchsian group, let $ \Gamma\vartriangleleft  \Gamma_{0}$   and let $\delta \geq  \delta(\Gamma)$. Then, the following mapping $\psi$ is a affine bijection between the Radon eigenmeasures of $L_{\phi^{X,\delta}}$ for eigenvalue $1$ and the $(\Gamma, \delta)$-conformal measures: For  a Radon eigenmeasure $\mu_X$ and a Borel set $E \subseteq \partial \mathbb{D}$, 
$$\psi(\mu_X)(E)  = \mu_X\bigl( \pi_{\Sigma}^{-1}(E) \times \{\Gamma\} \bigr). $$
\end{numberedtheorem}

The theory of the eigenmeasures of the Ruelle operator is well developed, see \cite{bowen_1975,sarig_1999,sarig_2001,mauldin_2001,stadlbauer_2017,shwartz_2019}. In particular, in \cite{shwartz_2019} the author presented the eigenmeasures of a transient Ruelle operator (see definition in Section \ref{sec:martin_boundary}) in terms of points on a Martin boundary. Thus, the classification of the conformal measures translates to the identification of the Martin boundary. 

\paragraph{Conformal measures and hyperbolic covers.}

In Sections \ref{sec:conf_of_hyperbolic_covers} and \ref{sec:conv_geo_flow} we apply the principle described above to the case where the group of deck transformations $G= \Gamma_0 / \Gamma$ is hyperbolic. In the canonical probabilistic setting, Ancona's well known theorem \cite{ancona_1987,ancona_1988} relates the Martin boundary of a finite range random walk  on a hyperbolic graph to the hyperbolic boundary of the graph. Using an extended version of Ancona's theorem to the Ruelle operator setting (see   Section \ref{sec:ancona_thm}) for every $\delta >\delta(\Gamma)$ we relate the $(\Gamma, \delta)$-conformal measures to the hyperbolic boundary of $G$, denoted by $\partial G$.     

In what follows,  a sequence $(a_i)$ with $a_i\in S_0$ is called a  \textit{boundary expansion }of a point $\xi\in \partial \mathbb{D}$ if for every $n\geq 0$, $f^n_{\Gamma_0}(\xi)\in I_{a_n}$. Observe that $(a_i)$ is a boundary expansion of a point $\xi\in \partial \mathbb{D}$ iff $\pi_{\Sigma}(a_0, a_1,\dots) = \xi$.  
\begin{numberedtheorem}[\ref{thm:conv_along_paths}]
Let $\Gamma_0$ be a cocompact Fuchsian group, let $ \Gamma\vartriangleleft  \Gamma_{0}$ and let $\delta > \delta(\Gamma)$. Assume that  $G=\Gamma_0 / \Gamma$ is a hyperbolic group.  Then, for every  $\mu \in \confG$,  for $\mu$-a.e. $\xi\in \partial \mathbb{D}$ with Bowen-Series coding $(a_n)$, the  sequence $$e_{a_{n}}^{-1}\dots e_{a_{0}}^{-1} \Gamma$$ converges to a point in $\partial G$. 
If $\mu\in \extconf$, then there exists $\eta\in \partial G$  s.t. the sequence almost-surely converges to $\eta$. Conversely, for every $\eta \in \partial G$, there exists a unique $\mu \in \extconf$ with $\eta$ its almost-surely limiting point of the sequence.
\end{numberedtheorem}

\begin{figure}[h!]
\centering
\includegraphics[width=2.5in,height=2.5in]{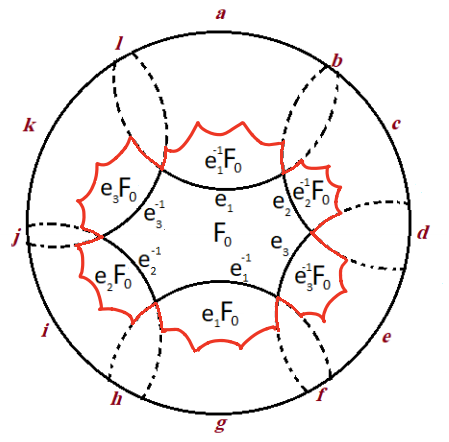}
\caption{In this example,  $S_{0} = \{a,b,c,d,e,f,g,h,i,j,k,l\}$, $e_a=e_{b} =e_1, e_c =e_{d}= e_2, e_e = e_{f}=e_3, e_g =e_{h}= e_1^{-1}, e_i = e_{j}=e_2^{-1}$ and $e_k = e_l=e_3^{-1}$. The elements $\{e_i\}$ map $F_0$ onto neighbouring copies. In particular, the copy $e_i^{-1} F_0$  shares with $F_0$ the edge labeled $e_i$. 
}
\label{fig:bs_coding_example}
\end{figure}

We derive a similar  result for cutting sequences of geodesics. In 
more details, let
$$\mathcal{R} = \bigl\{(\xi^-, \xi^+) \in (\partial\mathbb{D})^2: \text{ the geodesic curve between } \xi^- \text{ and }\xi^+ \text{ intersects }int F_0 \bigr\}. $$
 Since the group $\Gamma_0$ is cocompact, $F_0$ is a polygon
in $\mathbb{D}$
with finite number of edges.  For every $\gamma_1,\gamma_2\in \Gamma_0$,
$$int(\gamma_1 F_0) \cap int(\gamma_2 F_0) \neq \varnothing\Longleftrightarrow\
\gamma_1 = \gamma_2$$ and  $$\gamma_1 F_0 \text{ and }\gamma_2 F_0 \text{
share a common edge } \Longleftrightarrow\ \gamma_1 \gamma_2^{-1}\in \{e_a\}_{a\in
S_0}.$$See for example  Figure \ref{fig:bs_coding_example}. Given $(\xi^-, \xi^+)\in \mathcal{R}$, let $(F_{i})_{i\in \mathbb{Z}}$ be the sequence of copies of $F_0$ that the geodesic curve between $\xi^-$ and $\xi^+$ intersects. In case the curve passes through a vertex of some $F_i$, we perturb the curve around it, see Figure 5 in \cite{series_1986}.   Then, for all $i$ there exists a unique $e_i\in \{e_a\}_{a\in S_0}$ s.t. $ F_i = e_i^{-1} F_{i+1}$. The sequence $(e_i)$ is called the \textit{cutting sequence }of $(\xi^-, \xi^+)$. For $(\xi^-, \xi^+) \in \mathcal{R}$ with a cutting sequence $(e_{i})$, we write
$$\eta_n^+(\xi^-, \xi^+) :=e_{n}^{-1}\cdots e_0^{-1} \Gamma$$ 
and
$$\eta_{n}^-(\xi^-, \xi^+) :=(e_{-n-1})^{-1}\cdots (e_{-1})^{-1} \Gamma.$$
In Section  \ref{sec:conv_geo_flow}
 we prove the following theorem which describes the almost surely limiting behaviour of $\eta_n^+$ and $\eta_n^-$.

\begin{numberedtheorem}[\ref{thm:conv_cutting_seq}]
Let $\Gamma_0$ be a cocompact Fuchsian group, let $ \Gamma\vartriangleleft  \Gamma_{0}$, let $\delta > \delta(\Gamma)$ and let $\mu\in \confG$. Assume that $G=\Gamma/\Gamma_0$  is a hyperbolic group.   Then, 
\begin{enumerate}
\item 
For $ \mu$-a.e. $\xi^+\in\partial \mathbb{D}$, for every $\xi^-\in\mathbb{D}$ s.t. $(\xi^-, \xi^+)\in \mathcal{R}$, the sequence $\eta_{n}^+(\xi^-, \xi^+)$ converges to a point in $\partial G$.  
\item 
For $ \mu$-a.e. $\xi^-\in\partial \mathbb{D}$, for every $\xi^+\in\mathbb{D}$ s.t. $(\xi^-, \xi^+)\in \mathcal{R}$, the sequence $\eta_{n}^-(\xi^-, \xi^+)$ converges   to a point in  $\partial G$.  

\end{enumerate}
If $\mu\in \extconf$  then there exists a point $\eta\in \partial G$ s.t. the sequences  almost-surely converges to $\eta$. Conversely, for every $\eta \in \partial G$, there exists a unique $\mu \in \extconf$ with $\eta$ its almost-surely limiting point of the sequences. 
\end{numberedtheorem}

The limiting point $\eta$ from Theorem \ref{thm:conv_cutting_seq} is the same limiting point from Theorem \ref{thm:conv_along_paths}. We emphasize that Theorem \ref{thm:conv_cutting_seq} does not follow directly from Theorem \ref{thm:conv_along_paths} because when $\Gamma_0$ is cocompact the set of cutting sequences is \textbf{not} a Markov shift, see \cite{series_1986}. To prove the theorem, we use the description of cutting sequences developed by Series in \cite{series_1986}.

The extended Ancona's theorem (see  Theorem  \ref{thm:hyperbolic_main_thm})  holds only for supercritical values, which translates in this setting to $\delta >\delta(\Gamma)$. We were recently informed by Bispo and Stadlbauer that they can show that for a   potential function with a quasi-symmetric Green's function on a group extension of a hyperbolic group,  the results of the extended Ancona's theorem  also hold  at the critical value. Following this, we conjecture that  Theorems \ref{thm:conv_along_paths} and \ref{thm:conv_cutting_seq} should hold  at the critical value $\delta = \delta(\Gamma)$ as well, if the Poincar\'e series converges at the critical value.

For similar results on dependent random walks involving invariant measures  (rather than conformal measures), see  \cite{karlsson_1999,karlsson_2006}.  

\section{Preliminaries}
\label{sec:Preliminaries}

\subsection{Topological Markov shifts and the Ruelle operator}
Let $S$ be an infinite countable set of states and let $\mathbb{A} = (\mathbb{A}_{a,b})_{S\times S}\in \{0,1\}^{S\times S}$ be a transition matrix over $S$. For a subset  $A \subseteq \mathbb{Z}$ and a vector $x\in S^{A}$, we denote by $x_i$ the $i$-th coordinate of $x$. 

The \textit{(positive) one-sided topological Markov shift} (TMS) is the space
$$X^{} =\{x\in S^{\mathbb{N}\cup \{0\}} : \mathbb{A}_{x_i, x_{i+1} }= 1, \forall i\geq 0\} $$
with the transformation $T:X^{}\rightarrow X^{}$, $(Tx)_i = x_{i+1}$ and the  metric
\begin{equation*}
\label{eq:def_d}
d(x,y) = 2^{-\inf\{i\geq 0:x_i\neq y_i\}}.
\end{equation*}
If $\sum_{b}\mathbb{A}_{a,b}<\infty$ for every $a\in S$, then the space $(X^{},d)$ is locally compact  and all \textit{cylinder} sets
\begin{equation*}
\label{eq:def_cylinders}
[a_0,\dots, a_m] := \{x\in X^{} : x_i = a_i, 0\leq i \leq m\}
\end{equation*}
are compact. A word $(a_1, \dots, a_n)\in S^{n}$ is called \textit{admissible} if $[a_1,\dots, a_n]\neq \varnothing$. We denote by $\mathcal{W}^n$ the set of all admissible words of length $n$,
$$\mathcal{W}^{n}=\{(a_1,\dots, a_n)\in S^n : [a_1,\dots, a_n]\neq \varnothing\}. $$
We say that  $X^{}$ is \textit{topologically transitive}, or simply \textit{transitive}, if for every $a,b\in S$ there exists $n\geq 0$ s.t. $T^{-n}[a] \cap [b] \neq \varnothing$. 

For  numbers $r_1,r_2,c\in \mathbb{R}^+$, we write $r_1 =e^{ \pm c} r_2$ if
$ e^{-c}r_2 \leq r_1 \leq e^{c} r_2$. Given two functions $f,g:\mathbf{D}\rightarrow \mathbb{R}^+$ (or measures), where $\mathcal{\mathbf{D}}$ is some domain, we write $f \ll\ g$ if there exists $c>1$ s.t. $f (p)\leq c g(p)$ for all $p\in \mathcal{\mathbf{D}}$. We write $f \asymp g$ if $f \ll g$ and $g \ll f$.

We denote by $C_c(X^{})$  the space of all continuous functions from $X$ to $\mathbb{R}$ with compact support, by $C^+(X^{})$  the space of all non-negative continuous functions and by $C^+_c(X^{})=C^+(X) \cap C_c(X^{})$ the space of all non-negative continuous functions with compact support. 

The $m$-th \textit{variation} of a  function $\phi:X^{}\rightarrow\mathbb{R} $ is 
\begin{equation*}
\label{eq:def_var} Var_m(\phi) = \inf\{|\phi(x) - \phi(y)| : x,y\in X^{}, x_i=y_i, 0\leq i< m-1\}.
\end{equation*}
A function $\phi$  is said to have \textit{summable variations} if   $\sum_{m\geq 2} Var_m(\phi) <\infty$.
We let $\phi_n = \sum_{i=0}^{n-1}\phi \circ T^{i}$ and $C_\phi =\sum_{m\geq 2} Var_m (\phi)$.
\begin{defn}
\label{def:ruelle_operator}\em
 The \textit{Ruelle operator} $L_\phi$ evaluated on a function $f\in C(X^{})$ at a point $x\in X$ is
$$(L_\phi f)(x) = \sum_{y:Ty=x}e^{\phi(y)}f(y). $$
When $X$ is locally compact, the  sum is finite  for every $f\in C_{c}(X)$.
Then, for every $n>0$, 
$$(L_\phi^n f)(x)= \sum_{y:T^n y=x}e^{\phi_n(y)}f(y). $$
\end{defn}
\begin{defn}\em
The \textit{Gurevich pressure} of $\phi$ is the following limit
$$P_G(\phi) = \limsup _{n\rightarrow\infty}\frac{1}{n}\log \sum_{T^n x=x}e^{\phi_n(x)}1_{[a]}(x) $$
for some $a\in S$ and $x\in X^{}$. \end{defn}
If $(X,T)$ is topologically transitive and $\phi$ has summable variations, then $P_G(\phi)$ is independent of  the choice
of $a$, see \cite{sarig_1999}.
When $P_G(\phi)<\infty$, we write $\rho(\phi) = \exp\bigl( P_G(\phi)\bigr)$. 

\subsection{The Martin boundary of a transient potential.}
\label{sec:martin_boundary}
 
 Assume that $X$ is transitive and locally compact and that $\rho(\phi)< \infty$. Let $t\in \bigl[\rho(\phi) ,\infty\bigr)$. The $t$-\textit{Green's function}, evaluated at $f\in C_c(X)$ and $x\in X$, is 
$$G(f,x|t) := \sum_{n\geq0}t^{-n} (L_\phi^n f)(x).$$
We say that $\phi$ is $t$-\textit{recurrent} if $G(f,x|t)=\infty$ for some (or equivalently for every) $0 \not\equiv f\in C_c^+(X)$ and $x\in X$. Otherwise, we say that $\phi$ is $t$-\textit{transient}. If $\phi$ is $1$-transient, we simply say that it is \textit{transient}. Then, we write $G(f,x) := G(f,x|1)$. 
Notice that the ``transience'' in \cite{sarig_2001} means in our terminology $\rho(\phi)$-transience.

For a $t$-transient potential with summable variations, the author introduced in \cite{shwartz_2019} a Martin boundary which represents all eigenmeasures (and analogously all eigenfunctions) of the Ruelle operator, for the eigenvalue $t$. We briefly describe the construction here. 

Fix  $o\in S$ arbitrarily. When $\phi$ is $t$-transient, for every fixed $f\in C_c^+(X)$, the \textit{Martin kernel } $$K(f,x|t) := \frac{G(f,x|t)}{G(1_{[o]},x|t)}$$ is  continuous and bounded as a function of $x$. Let $\{w_i\}_{i\in \mathbb{N}}$ be an enumeration of $\cup_{i\geq 1}\mathcal{W}_i$. We define a new metric on $X$, 
$$\varrho(x,y|t) = \sum_{i=1 }^\infty \frac{|K(1_{[w_i]},x|t) -K(1_{[w_i]},y|t)| + |1_{[w_{i}]}(x) -1_{[w_{i}]}(y)|}{\max_{z\in [w_i]}|K(1_{[w_i]},z|t) |}. $$
The \textit{$t$-Martin compactification}, denoted by $\widehat X(t)$, is the completion of  $X$ w.r.t. the metric $\varrho$. The \textit{$t$-Martin boundary}, denoted by $\mathcal{M}(t)$, is the set of  all  new obtained points, $\mathcal{M}(t) = \widehat X(t)\setminus X$. For every fixed $f\in C_c(X)$, the Martin kernel $K(f,\cdot|t)$ is a $\varrho$-continuous function in $X$ and it can be uniquely extended to a $\varrho$-continuous function in $\widehat{X}(t)$ via
$$K(f,\omega|t) = \lim_{x\rightarrow\omega}K(f,x|t), \quad \omega \in \mathcal{M}(t). $$ 
For $\omega \in \mathcal{M}(t)$ and $f\in C_c(X)$, we let $\mu_\omega(f) := K(f, \omega|t)$. Observe that for every $\omega \in \mathcal{M}(t)$, the measure $\mu_\omega$ is a $t$-eigenmeasure of $L_\phi$.

The \textit{$t$-minimal boundary}  $\mathcal{M}_m(t)$ is the set of all points  $\omega \in \mathcal{M}(t)$ s.t. the resulting measure $\mu_\omega$ is extremal in the cone of  eigenmeasures for eigenvalue $t$.
Then, for every positive Radon measure $\mu$ with $L_\phi^*\mu = t\mu$ there exists a unique finite measure $\nu$ on $\mathcal{M}_m(t)$ s.t.
\begin{equation}
\label{eq:minimal_ker_rep}
\mu(f) = \int_{\mathcal{M}_m(t)} \mu_\omega (f)d\nu(\omega), \quad \forall f\in C_c(X).
\end{equation}

By definition, a sequence $x^n\in X$ converges to a point $\omega\in \mathcal{M}(t)$ in the topology of $\widehat X(t)$ iff $K(f,x^n|t) \xrightarrow[n\rightarrow\infty]{} K(f,\omega|t)$ for all $f\in C_c(X)$. In particular, 
a  point $\omega \in \mathcal{M}_m(t)$ is fully characterized by the following convergence property: for $\mu_\omega$-a.e. $x\in X$, $T^nx\rightarrow \omega$ in $\widehat{X}(t)$.

 In this paper, we mainly assume that $P_G(\phi) < 0$, which directly implies that $\phi$ is transient and the Martin boundary $\mathcal{M} :=\mathcal{M}(1)$ exists.
We write $\mathcal{M}_m:=\mathcal{M}_m(1)$.

\subsection{The hyperbolic boundary}
\label{sec:hyperbolic_boundary}
We briefly recall  the definitions of a hyperbolic graph and its boundary. For more detailed description, see  \cite{ghys_1990} and also \cite{woess_2000}.
\begin{defn}\em
Let $E \subseteq S\times S$ be a set of edges over $S$. We say that $(S,E)$ is \textit{connected} if for every $a,b\in S$, there exist $a_1, \dots, a_n\in S$ s.t. $a_1=a, a_n=b$ and  $(a_i, a_{i+1})\in E$, $1\leq i <n$. We say that $(S,E)$ is \textit{undirected}  if $$(a,b)\in E\Longleftrightarrow(b,a)\in E.$$
We say that $E$ is  \textit{locally finite}  if for every $a\in S$, $\#\{b:(a,b)\in E\}<\infty$.
 We denote by $d_E(a,b)$ the length of a shortest path from $a$ to $b$ in $(S,E)$. When $(S,E)$ is undirected and connected, $d_E$ is a metric.
\end{defn}
\begin{defn}\em
Let $(S,E)$ be a connected, undirected and locally finite graph. A \textit{geodesic triangle} consists of three points $a,b,c\in S$ and three geodesic paths  $\pi(a,b), \pi(b,c), \pi(c,a)$ from $a$ to $b$, $b$ to $c$ and $c$ to $a$ respectively.   
We say that the graph $(S,E)$ is  $\delta$\textit{-hyperbolic }if every geodesic triangle in the graph is \textit{$\delta$-thin}, namely any point on one of its sides is at distance at most $\delta$ from the other two sides.
\end{defn}
Let $(S,E)$ be a $\delta$-hyperbolic graph. For $a,b,o\in S$, let
$$|a\wedge b| _{o}:= \frac{1}{2}\bigl(d_{E}(o,a)+ d_{E}(o,b)- d_{E}(a,b)\bigr).
 $$
Fix an origin point $o\in S$. 
\begin{defn}\em
We say that a sequence $a_n\in S$ \textit{converges to the hyperbolic boundary} in $(S,E)$ if $$\lim_{m,n\rightarrow\infty}|a_n\wedge a_m|_{o}=\infty.  $$
Two sequences converging to the hyperbolic boundary $a_n,b_n\in S$ are said to be equivalent if 
\begin{equation}
\label{eq:hyperbolic_relation}
\lim_{n\rightarrow\infty}|a_n \wedge b_n|_{o} \rightarrow\infty. \end{equation}
Easy to verify that these definitions do not depend on $o\in S$. 
\end{defn}
 
\begin{defn}\em
The \textit{hyperbolic boundary} (or the \textit{Gromov boundary}) of $(S,E)$, denoted by $\partial (S,E)$, is the collection of all equivalence classes according to the relation in  Eq. (\ref{eq:hyperbolic_relation}). \end{defn}

\subsection{Generalized Ancona's theorem}
\label{sec:ancona_thm}

Recall that  $\rho(\phi)$ is the radius of convergence of the Green's function. Consider a finite range random walk on a hyperbolic graph.  Then, for all $\lambda > \rho(\phi)$ the minimal $\lambda$-Martin boundary of the walk and the hyperbolic boundary of the graph coincide. It was first proven by Series for random walks on Fuchsian groups  \cite{series_1983}
and then by Ancona for more general hyperbolic graphs \cite{ancona_1987,ancona_1988}. See also \cite{kaimanovich_1994} for a similar result in more general spaces and  \cite{blachere_2011} for the connection between Ancona's inequality and the Green metric.
Later on, this result was   proved at the critical value $\lambda=\rho(\phi)$  by Gou\"ezel and  Lalley for random walks on Fuchsain groups \cite{gouezel_2012} and by  Gou\"ezel for  symmetric random walks on hyperbolic groups   \cite{gouezel_2014}.

To prove the main results of Sections \ref{sec:conf_of_hyperbolic_covers} and \ref{sec:conv_geo_flow}, we extend Ancona's theorem to the Ruelle operator setting.
In particular, we show that if  the potential is uniformly irreducible w.r.t. a hyperbolic graph (see Definition \ref{def:uniform_irr}) then for all $\lambda > \rho(\phi)$, the minimal Martin boundary $\mathcal{M}_m(\lambda)$ and the hyperbolic boundary coincide. The proof of the theorem, which is of technical flavour, appears in the appendix.

\begin{defn}\em
\label{def:uniform_irr}
Let $(S,E)$ be a connected, undirected and locally finite graph.
 We say that  $\phi$ is  \textit{uniformly irreducible} w.r.t. $(S,E)$ if:

\begin{enumerate}
\item 
$\phi$ is bounded;\item
For every $a,b\in S$ with $[a,b]\neq \varnothing$ we have that $(a,b)\in E;$ 
\item
  There exists $K>0$ s.t. for every $a,b\in S$ with $(a,b)\in E$, there exists $k\leq K$ with 
$$L_\phi^k1_{[a]}(bx_b)> 0.$$
\end{enumerate}
\end{defn}
\begin{remark} \em{}There may be $(a,b)\in E$ with $[a,b]=\varnothing$. The set of edges $E$ is symmetric and  we may have different values of $k$ for $(a,b)$ and $(b,a)$. However,  both values are still bounded by $K$.  \end{remark}

For every $a\in S$, we pick $x_{a}\in T[a]$ arbitrarily.

\begin{thm}
\label{thm:hyperbolic_main_thm}
Assume the following:
\begin{itemize}
\item 
$(X,T)$ is locally compact and topologically transitive.
\item
$\phi$ has summable variations and  $P_G(\phi)<\infty$. 
\item 
There exist $\delta \geq 0$ and a   $\delta$-hyperbolic graph $(S,E)$ s.t. 
$\phi$ is uniformly irreducible w.r.t.   $(S,E)$.
\end{itemize}
Then,  for every $\lambda > \rho(\phi)$, there is a  bijection  $\omega:\partial (S,E) \rightarrow \mathcal{M}_m(\lambda)$ s.t. for any $a_n\in S$,  $$a_n\xrightarrow[n\rightarrow\infty]{} \xi \in \partial (S,E)  \quad\Longleftrightarrow\quad \forall f\in C_{c}(X), \;K(f,a_nx_{a_n}|\lambda)\xrightarrow[n\rightarrow\infty]{} K(f,\omega(\xi)|\lambda).
$$
\end{thm}
We were recently informed by Bispo and Stadlbauer that they can show that if $X$ is a \textbf{group extension }of a hyperbolic group and the Green's function is quasi-symmetric  then the result of Theorem  \ref{thm:hyperbolic_main_thm} also holds at the critical value $\lambda = \rho(\phi)$.

%% file: hyperbolic.tex
\subsection{Regular covers of compact hyperbolic surfaces}
\label{sec:HyperbolicCovers}

Recall that $\mathbb{D} = \{z\in \mathbb{C}:|z|<1\}$ is the unit open hyperbolic disc and that $\partial \mathbb{D} = \{z\in \mathbb{C}: |z|=1\}$ is its boundary. We denote by $d_{\mathbb{D}}$  the hyperbolic metric on $\mathbb{D}$.
A Fuchsian group $\Gamma_{0}$ is said to be \textit{co-compact} if  $\mathbb{D}/
\Gamma_0$ is a compact surface. A \textit{regular cover} of $\mathbb{D} / \Gamma_0$ is a surface $\mathbb{D}/\Gamma$ where $\Gamma \lhd \Gamma_0$. The group of deck transformations $G$ can be identified with $ \Gamma_0/\Gamma$ as follows: $\gamma \Gamma \cdot x\Gamma = \gamma x \Gamma$, with $\gamma\in \Gamma_0$ and $x\in \mathbb{D}$.  Let $Fix(\Gamma_0) = \bigl\{\xi \in \partial \mathbb{D}: \exists \gamma_0\in \Gamma_0\setminus\{id\} \text{ s.t. }\gamma_0 \xi = \xi\bigr \}$. Notice that $Fix(\Gamma_0)$ is a countable set.

We denote by $\delta(\Gamma)$  the \textit{critical exponent} of $\Gamma$, namely the critical value of $\delta$ s.t. the \textit{Poincar\'e series} \begin{equation}
\label{eq:poincare_series}
  \mathbf{p}(\Gamma, \delta) :=\ \sum_{\gamma\in \Gamma}e^{-\delta d_{\mathbb{D}}(0, \gamma 0)} 
\end{equation}
converges for all $\delta > \delta(\Gamma)$ and diverges for all $\delta < \delta (\Gamma)$. In general, $\delta(\Gamma) \leq \delta(\Gamma_0)$ and there is an equality iff $G$ is amenable \cite{stadlbauer_2013}. See also \cite{jaerisch_2015,dougall_2016,coulon_2018} for similar results in more general spaces. Since $\Gamma_0$ is co-compact, $\delta(\Gamma_0)=1$, see Theorem 1.6.3 in \cite{nicholls_1989}. 

 We remind the reader the definition of a conformal measure:
\begin{defn}[Sullivan \cite{sullivan_1979}]\em
A finite positive measure $\mu$ on $\partial \mathbb{D}$ is said to be \textit{$(\Gamma, \delta)$-conformal} if for every $\gamma \in \Gamma$, 
$$ \frac{d(\mu \circ \gamma)}{d\mu} = |\gamma'|^\delta $$
where $(\mu \circ \gamma )(A) = \mu(\gamma A) = \int1_{A}(\gamma^{-1}x)d\mu(x)$.
We denote by $\confG$ the collection of all  $(\Gamma, \delta)$-conformal measures and by $\extconf$ the extremal points of $\confG$.
\end{defn}
Patterson and Sullivan originally considered what they called ``\textit{conformal densities}'' rather than conformal measures. However, both  definitions coincide, see Remark 3.3 in \cite{babillot_2004}.

\subsection{The Bowen-Series coding and its group extension}
\label{sec:bowen_series_coding}
For a  cocompact Fuchsian group $\Gamma_0$ with a fundamental domain $F_0\subseteq\mathbb{D}$, Bowen and Series constructed (w.r.t. $F_0$) in \cite{bowen_1979} a finite partition $\{I_a\}_{a\in S_0}$ of $\partial \mathbb{D}$ into closed arcs with disjoint interiors, a  finite set $\{e_a\}_{a\in S_0}\subseteq\Gamma_0$ and a map $f_{\Gamma_0}:\partial \mathbb{D}\rightarrow \partial\mathbb{D}$ with the following properties:
\begin{description}
\item[(Gen)]
The set $\{e_a\}_{a\in S_0}$ is symmetric and generates $\Gamma_0$. 
\item[(Res)]
For all $a\in S_0$, 
$f_{\Gamma_0}= e_a^{-1}$ on $in t(I_a)$. 
\item[(Mar)]
$\{I_a\}$ is a \textit{Markov partition}: if $int(f_{\Gamma_0}(I_a)) \cap int(I_b)\neq \varnothing $ then $I_b \subseteq f(I_a)$.
\item[(Tr)]
For every $a,b\in S_0$ there exists $n$ s.t. $f_{\Gamma_0}^n (I_a) \supseteq I_b$. 
         
\item[(Orb)]
For all except finitely many $\xi_1, \xi_2\in \partial \mathbb{D}$: 
$$\exists  n,m\in \mathbb{N} \text{ s.t. }f_{\Gamma_0}^n (\xi_1) = f^m_{\Gamma_0}(\xi_2) \quad \Longleftrightarrow\quad\exists \gamma_0\in \Gamma_0 \text{ s.t. }\xi_1 = \gamma_0 (\xi_2).$$

\item[(Bnd)]
There is a constant $N_0$ s.t. for every $\gamma_0\in \Gamma_0$, $\gamma_0\neq id_{\Gamma_0}$,
$$\# \bigcup_{n\geq 0}\left\{ (a_0,\dots, a_n)\in S_{0}^{n+1} : int(\cap_{i=0}^n f_{\Gamma_0}^{-i}(I_{a_i}))\neq \varnothing \text{ and }\gamma_0 = e_{a_n}^{-1}\cdots e_{a_0}^{-1} \right\} \leq N_0. $$

\item[(Dist)]
There exists a constant $B>1$ s.t. for every $a_1,\dots, a_n \in S_{0}$ and every $\xi_1,\xi_2\in \partial \mathbb{D}$ with $f_{\Gamma_0}^k \xi_i = e_{a_k}^{-1}\dots e_{a_1}^{-1}\xi_i$ for $k=1,\dots ,n$ and $i=1,2$, 
$$\frac{|(f_{\Gamma_0}^{n-1})'(\xi_1)|}{|(f_{\Gamma_0}^{n-1})'(\xi_2)|}\leq B .$$

\end{description}
For \textbf{(Orb)} and \textbf{(Bnd)}, see \cite{series_1981}.
For \textbf{(Dist)} see also \cite{ledrappier_2007}.

We write $$I_{a_1,\dots, a_n} = \cap_{i=1}^{n}f_{\Gamma_0}^{-i+1}I_{a_i} .$$ In particular, for all $\xi\in int( I_{a_1,\dots, a_n})$ and $1\leq k\leq n$, $$f_{\Gamma_0}^k(\xi )=  e_{a_k}^{-1}\cdots e^{-1}_{a_1}\xi.$$
For an admissible word $w=(a_1,\dots, a_n$), we write
$$e_{w} = e_{a_1}\dots e_{a_n}. $$

\begin{defn} \em
A sequence $(a_i)_{i\geq 0}$ with $a_i\in S_{0}$ is called a \textit{boundary expansion} of a point $\xi\in \partial \mathbb{D}$ if for every $n\geq 0$, $f_{\Gamma_0}^{n}(\xi) \in I_{a_n}$. 
\end{defn}
Let
\begin{equation}
\label{eq:pos_sidede_shift}
\Sigma=\{(\sigma_i)_{i\geq 0}  : \forall i\geq 0,\; \sigma_i \in S_{0}\text{ and } int (I_{\sigma_{i+1}}) \subseteq  int( f(I_{\sigma_i})) \} 
\end{equation}
and let $T_\Sigma:\Sigma \rightarrow\Sigma$ be the left-shift. 
Let $\pi_{\Sigma}:\Sigma\rightarrow\partial\mathbb{D}$  be the canonical projection,  $\pi_\Sigma(\sigma) \in \cap_{n\geq 0}\overline{f_{\Gamma_0}^{-n}I_{\sigma_n}}$ (the intersection is a singleton, see \cite{bedford_1991}). By \textbf{(Mar)} and \textbf{(Tr)}, $(\Sigma, T)$ is a one-sided transitive TMS and by \textbf{(Res)}, $f_{\Gamma_0}\circ \pi_\Sigma = \pi_\Sigma \circ T_\Sigma$.
Given a point $\sigma\in \Sigma$ we write $\sigma_i$ for its $i$-th coordinate.

Let $G = \Gamma_0 / \Gamma$ where $\Gamma \vartriangleleft \Gamma_0$. Let $(X,T)$ be the one-sided TMS over the set of states $S_{X}=S_{0}\times G$ with the following transition rule
\begin{equation}
\label{eq:X_transition}
(a,\gamma_1\Gamma) \rightsquigarrow (b, \gamma_2\Gamma) \Longleftrightarrow int(I_{{b}}) \subseteq  int (f_{\Gamma_0}(I_a))  \text{ and }\gamma_2\Gamma = e_a ^{-1}\gamma_1\Gamma. 
\end{equation}
The shift space $(X,T)$ is called  the \textit{group extension}, or the \textit{$G$-extension},  of $\Sigma$, see \cite{stadlbauer_2013}. We denote by $\pi_X:X\rightarrow\Sigma$ the natural projection from $X$ to $\Sigma$. 
\begin{defn}\em
Given $\delta>0$, let $\phi^{
\Sigma,\delta}:\Sigma\rightarrow\mathbb{R}$,  
$$\phi^{\Sigma,\delta}(\sigma) := -\delta_{}\log|(e_{\sigma_0}^{-1})'(\pi_{\Sigma}(\sigma))| $$ 
and let
$$\phi^{X,\delta}(x): = \phi^{\Sigma,\delta}(\pi_X(x)).$$
\end{defn}
  \begin{prop}[Series \cite{series_1981}] $\phi^{\Sigma,\delta}$ is H\"older continuous.
\end{prop}
Clearly $\phi^{X,\delta}$ is H\"older continuous as well.

The following propositions are elementary. For completeness, we provide their proofs in the appendix.

\begin{prop}
\label{prop:group_ext_transitive}
$(X,T)$ is topologically transitive.
\end{prop}

\begin{prop}
\label{prop:pressure_delta}
\begin{enumerate}
\item 
The potential $\phi^{X,\delta}$ is transient iff $\mathbf{p}(\Gamma, \delta)<\infty$.
\item  For every $\delta > \delta(\Gamma)$, 
$$P_G(\phi^{X, \delta})<0 .$$
\end{enumerate} \end{prop}

The following proposition allows us to exclude from our discussion measures with atoms. Its proof included in the appendix as well.
אד פר\begin{prop}
\label{prop:non_atomic}
Let $\delta \geq \delta(\Gamma)$. 
\begin{enumerate}[label={(\arabic*)}]
\item 
Every   $(\Gamma, \delta)$-conformal measure is non-atomic.
\item
Every Radon eigenmeasure of $L_{\phi^{X,\delta}}$ for eigenvalue $1$ is non-atomic. 
\end{enumerate}
\end{prop}
\hfill\ensuremath{\square}

\section{The eigenmeasures of the Ruelle operator and  the conformal measures}
\label{sec:conf_bijection}
In this section we relate the $(\Gamma,\delta)$-conformal measures to the eigenmeasures of $L_{\phi^{X,\delta}}$ for eigenvalue $1$.   

\begin{thm}
\label{prop:gamma_conf_bijection}
Let $\Gamma_0$ be a cocompact Fuchsian group, let $ \Gamma\vartriangleleft  \Gamma_{0}$   and let $\delta >  \delta(\Gamma)$. Then, the following mapping $\psi$ is a affine bijection between the Radon eigenmeasures of $L_{\phi^{X,\delta}}$ for eigenvalue $1$ and the $(\Gamma, \delta)$-conformal measures: For  a Radon eigenmeasure $\mu_X$ and a Borel set $E \subseteq \partial \mathbb{D}$,
$$\psi(\mu_X)(E)  = \mu_X\bigl( \pi_{\Sigma}^{-1}(E) \times \{\Gamma\} \bigr). $$
\end{thm}

Before proving the theorem, we deduce several elementary results, some already known, using Theorem \ref{prop:gamma_conf_bijection} and known theory on the eigenmeasures of the Ruelle operator.

\begin{cor}
\label{cor:unique_integral_rep}
Let $\delta \geq\delta(\Gamma)$. Then, for every $\mu \in \confG$ there exists a unique finite measure $\nu$ on $\extconf$ s.t. 
$$\mu = \int_{\mu'\in \extconf}\mu' d\nu(\mu'). $$ 
\end{cor}
\begin{proof}
This result can be derived  from the classical Choquet theory but also  follows by the  unique representation on the minimal boundary (see Equation (\ref{eq:minimal_ker_rep}) and also  \cite{shwartz_2019})  and by the linearity of the mapping in Theorem \ref{prop:gamma_conf_bijection}. 
\end{proof}

\begin{cor}[Furstenberg \cite{furstenberg_1973}]
Assume that $\mathbb{D}/\Gamma$ is compact. Then,   a $(\Gamma, \delta)$-conformal measure exists iff $\delta = \delta(\Gamma)$. Moreover, the $(\Gamma, \delta(\Gamma))$-conformal measure is unique up to scaling.
\begin{proof}
The corollary follows directly from the Ruelle's Perron-Frobenius  theorem, see \cite{bowen_1975}.
\end{proof}
\end{cor}

\begin{cor}[C.f. Sullivan \cite{sullivan_1987}]
Assume that $\mathbf{p}(\Gamma, \delta(\Gamma))=\infty$. Then, the $(\Gamma, \delta(\Gamma))$-conformal measure is unique up to scaling.
\begin{proof}
By Proposition \ref{prop:pressure_delta}, the potential $\phi^{X, \delta(\Gamma)}$ is recurrent. By Sarig's generalized  Ruelle's Perron-Frobenius Theorem \cite{sarig_1999,sarig_2001}, the  eigenmeasure of $L_\phi^{X^\delta}$ is unique up to normalization. 
\end{proof}
\end{cor}

In what follows, let $Y = \partial\mathbb{D}\times G$. The group $\Gamma_0$ acts on $Y$ in the following way: 
$$\gamma_0(\xi, \gamma\Gamma) = (\gamma_0 \xi, \gamma_0 \gamma \Gamma), \quad \gamma_0 \in \Gamma_0, (\xi, \gamma\Gamma)\in Y.$$
Let
$ \ f_{Y}: Y\rightarrow  Y$ be the extension of $f_{\Gamma_0}$ to $Y$:$$ f_{Y}(\xi, \gamma\Gamma) = (e_a^{-1}\xi, e_a^{-1}\gamma\Gamma), \quad \xi \in int(I_a). $$
Since we narrowed our discussion to non-atomic measures, we may ignore the values of $f_Y$ on $\partial I_a$. To prove Theorem \ref{prop:gamma_conf_bijection}, we map, in several steps the Radon eigenmeasures of $L_{\phi}^{X,\delta}$ for eigenvalue $1$ to the Radon measures on $Y$ which satisfies a $\Gamma_0$-regularity condition, see Eq. (\ref{eq:mu_Y_gamma_prop}) in the following lemma. 
\begin{lemma}
\label{lemma:mu_Y_prop}
Let $\mu_Y$ be a non-atomic  Radon measure on $Y$. Then, the following are equivalent:\begin{enumerate}
\item 
The measure $\mu_Y\circ f_{\Gamma_0}$ given by  $$(\mu_Y \circ f_Y)(A\times \{\gamma\Gamma\}) = \sum_{a\in S_0 }\mu_Y\biggl(f_Y\bigl((I_a\cap A) \times \{\gamma\Gamma\}\bigr)\biggr)$$ with $A \subseteq \partial \mathbb{D}$ measurable is absolutely continuous w.r.t. $\mu_Y$ and
\begin{equation}
\label{eq:RN_f_Y}
\frac{d(\mu_Y \circ f_Y)}{d\mu_Y}(\xi, \gamma\Gamma) =   |f_{\Gamma_0}'(\xi) |^{\delta},   \quad \mu_Y -a.e. 
\end{equation}

\item
The measure $\mu_Y$ is $\Gamma_0$-quasi-invariant and for all $\gamma_0\in \Gamma_0$, 
\begin{equation}
\label{eq:mu_Y_gamma_prop}
\frac{d(\mu_Y \circ \gamma_0)}{d\mu_Y}(\xi, \gamma\Gamma) = |\gamma_0'(\xi) |^{\delta}, \quad \mu_Y-a.e.
 \end{equation}
\end{enumerate} 
\end{lemma}
\begin{proof} Assume (1) holds. 
Fix $\gamma_0\in \Gamma_0$, $\gamma_0 \neq id_{\Gamma_0}$ and let
$$A_{n,m} = \{\xi\in \partial \mathbb{D} : f_{\Gamma_0}^n(\xi) = ( f_{\Gamma_0}^m \circ \gamma_0)(\xi) \}. $$
 By \textbf{(Orb)},  $\mu_Y\bigl((\bigcup_{n,m\geq 0} A_{n,m} \times G)\bigtriangleup\ Y\bigr)=0$.   
 Fix $n,m\geq 0$ and let $\xi \in A_{n,m}\setminus Fix(\Gamma_0)$. Let $a_1,\dots, a_n, b_1,\dots, b_m\in S_{0}$ s.t. $\xi \in I_{a_1,\dots, a_n}$ and $\gamma_0 \xi \in I_{b_1,\dots, b_m}$. Then, 
$$e_{a_n}^{-1}\cdots e_{a_1}^{-1} (\xi)= e_{b_m}^{-1}\cdots e_{b_1}^{-1}\gamma_0 (\xi). $$ 
In particular, $\xi$ is a fixed point of $e_{a_1}\cdots e_{a_n}e_{b_m}^{-1}\cdots e_{b_1}^{-1} \gamma_0$. Since $\xi \not\in Fix(\Gamma_0)$, 
$$ \gamma_0 = \left( e_{b_m}^{-1}\cdots e_{b_1}^{-1} \right)^{-1}e_{a_n}^{-1} \cdots e_{a_1}^{-1}$$
and 
\begin{align*} |\gamma_0'(\xi)|^\delta =&\left|\left(  \left(e_{b_m}^{-1}\cdots e^{-1}_{b_1}\right)^{-1}\right)'(e_{a_n}^{-1}\cdots e_{a_2}^{-1}\xi)\right|^\delta\cdot|(e_{a_n}^{-1}\cdots e_{a_{1}}^{-1})'(\xi)|^\delta\\
=&
 \frac{|(e_{a_n}^{-1}\cdots e^{-1}_{a_1})'(\xi)|^{\delta}}{|(e_{b_m}^{-1}\cdots e^{-1}_{b_1})' (\gamma_0 \xi)|^\delta}. 
\end{align*}
For $\mu_Y$-a.e. $(\xi, \gamma\Gamma)\in  \left(I_{a_1,\dots, a_n}\cap \gamma_0^{-1} I_{b_1,\dots, b_m}\right)\times G$ we have that
\begin{align*}\frac{d(\mu_Y \circ f_Y^n)}{d\mu_Y}(\xi, \gamma\Gamma) = & \frac{d(\mu_Y \circ f_Y^m \circ\gamma_0) }{d\mu_Y} (\xi, \gamma\Gamma) \\
=& \frac{d(\mu_Y \circ f_Y^m)}{d\mu_Y}(\gamma_0\xi, \gamma_0 \gamma\Gamma)\frac{d(\mu_Y\circ \gamma_0)}{d\mu_Y}(\xi, \gamma\Gamma). 
\end{align*}
Moreover, 
 by Eq. (\ref{eq:RN_f_Y}), for $\mu_Y$-a.e. $(\xi, \gamma\Gamma) \in \left(I_{a_1,\dots, a_n}\cap \gamma_0^{-1} I_{b_1,\dots, b_m}\right)\times G$,
$$\frac{d(\mu_Y \circ f_Y^n)}{d\mu_Y}(\xi, \gamma\Gamma) = |(e_{a_n}^{-1}\cdots e_{a_{1}}^{-1})'(\xi)|^\delta $$
and
$$\frac{d(\mu_Y \circ f_Y^m)}{d\mu_Y}(\gamma_0\xi, \gamma_0 \gamma\Gamma) = |(e_{b_m}^{-1}\cdots e_{b_{1}}^{-1})'(\gamma_0\xi)|^\delta. $$
Thus, for $\mu_Y$-a.e. $(\xi, \gamma\Gamma) \in \left(I_{a_1,\dots, a_n}\cap \gamma_0^{-1} I_{b_1,\dots, b_m}\right)\times G$,
$$\frac{d(\mu_Y\circ \gamma_0)}{d\mu_Y}(\xi, \gamma\Gamma) =\left(\frac{d(\mu_Y \circ f_Y^m)}{d\mu_Y}(\gamma_0\xi, \gamma_0 \gamma\Gamma)\right)^{-1} \frac{d(\mu_Y \circ f_Y^n)}{d\mu_Y}(\xi, \gamma\Gamma) = |\gamma_0'(\xi)|^\delta. $$
Since there is only a countable number of such $(a_1,\dots, a_n), (b_1,\dots, b_m),m$ and $n $, the identity holds for $\mu_Y$-a.e. $(\xi, \gamma\Gamma)\in Y$.
So $(1)\Rightarrow (2)$.

Next, assume (2). Fix $a\in S_{0}$. Then, for $\mu_Y$-a.e. $(\xi, \gamma\Gamma)\in I_a \times G$,
$$\frac{d(\mu_Y \circ f_Y)}{d\mu_Y}(\xi, \gamma\Gamma) = \frac{d(\mu_Y\circ e_{a^{}}^{-1})}{d\mu_Y}(\xi, \gamma\Gamma) =|(e_{a}^{-1})'(\xi)|^\delta =|f_{\Gamma_0}'(\xi)|^\delta.  $$
\end{proof}

Henceforth we  use the following canonical correspondence to identify $X$ with $\Sigma \times G$,
$$(\sigma, \gamma\Gamma) \longmapsto \bigl((\sigma_0, \gamma\Gamma),(\sigma_1, e_{\sigma_0}^{-1}\gamma\Gamma),(\sigma_1, e_{\sigma_1}^{-1}\gamma\Gamma),\dots\bigr).  $$
In particular, we will not distinguish between the two. We let $\tilde \pi:X\rightarrow Y$, $\tilde \pi(\sigma, \gamma \Gamma) = (\pi_{\Sigma}(\sigma), \gamma \Gamma)$.

\begin{lemma}
\label{lemma:mu_X_mu_Y_bijection}
Let $\delta\geq \delta(\Gamma)$. Then, the map $\mu_X \mapsto  \mu_X \circ \tilde\pi^{-1}$ is an affine bijection between the Radon eigenmeasures of $L_{\phi^{X,\delta}}$ with eigenvalue $1$ and the  non-atomic Radon measures on $Y$ which satisfy Eq. (\ref{eq:mu_Y_gamma_prop}).
\end{lemma}
\begin{proof}
Recall that $\pi_\Sigma$ is bijective away from a countable number of points, see \cite{series_1981}. Therefore, since all eigenmeasures of the Ruelle operator are non-atomic (see Proposition \ref{prop:non_atomic}), $\tilde \pi$ is a measure-theoretic isomorphism.  

 Recall that $\mu_X$ is an eigenmeasure of $L_{\phi^{X, \delta}}$ of eigenvalue $1$ iff 
$$\frac{d(\mu_X \circ T)}{d\mu_X}(\sigma, \gamma\Gamma) =   |f_{\Gamma_0}'(\pi(\sigma))|^{\delta} $$
where ($\mu_X \circ T)(A\times \{\gamma\Gamma\}) =  \sum_{a\in S_{0}}\mu_X(T(([a]\cap A)\times \{\gamma \Gamma\})) $, see \cite{sarig_2015} and references within. Since $\tilde \pi \circ f_Y = T \circ \tilde \pi$, 
\begin{align*}\frac{d(\mu_Y \circ f_Y)}{d\mu_Y}\bigl( (\pi^{-1}(\sigma), \gamma\Gamma)\bigr )=& \frac{d(\mu_Y \circ (f_Y \circ \tilde \pi^{-1}))}{d(\mu_Y \circ \tilde \pi^{-1})}( \sigma, \gamma\Gamma)\\
 =
 & \frac{d(\mu_X \circ T)}{d\mu_X}(\sigma , \gamma\Gamma)\\
 =&|f_{\Gamma_0}'(\pi(\sigma))|^{\delta}  
\end{align*}
Hence, by Lemma \ref{lemma:mu_Y_prop},  $\mu_X$ is an eigenmeasure iff $\mu_Y$ satisfies Eq. (\ref{eq:mu_Y_gamma_prop}).
\end{proof}

\paragraph{\bf{Proof of Theorem \ref{prop:gamma_conf_bijection}.}}
By Lemma \ref{lemma:mu_X_mu_Y_bijection}, it suffices to present a bijection between the  $(\Gamma, \delta)$-conformal measures and the non-atomic Radon measures on $Y$ which satisfy Eq. (\ref{eq:mu_Y_gamma_prop}).

Let $\mu\in \confG$.
 We define a new measure $\mu_Y$ on $Y=\partial \mathbb{D} \times G$ as follows: For $A\subseteq \partial\mathbb{D}$ and $\gamma\Gamma\in G$,  
\begin{equation}
\label{eq:mu_Y_def}
\mu_Y(A\times \{\gamma\Gamma\}) :=\int |\gamma'(\xi)|^\delta 1_{A}(\gamma\xi)d\mu(\xi).
\end{equation} 
 We show that this  definition does not depend on the choice of $\gamma$ which represents $\gamma\Gamma$. Assume that $\gamma_1 \Gamma = \gamma_2 \Gamma$  
 and let $\gamma\in \Gamma$ s.t. $\gamma_1 = \gamma_2 \gamma$. Since $\mu$ is $(\Gamma,\delta)$-conformal,
\begin{align*}
 \int|\gamma_1
'(\xi)|^\delta \ 1_{A}(\gamma_1\xi)d\mu(\xi) = & \int  |(\gamma_2\gamma)
'(\xi)|^\delta 1_{A}(\gamma_2\gamma\xi)d\mu(\xi)\\
=&\int  |\gamma_2
'(\gamma\xi)|^\delta|\gamma'(\xi)|^\delta1_{A}(\gamma_2\gamma\xi)d\mu(\xi)\\
=& \int  |\gamma_2
'(\xi)|^\delta|\gamma'(\gamma^{-1}\xi)|^\delta1_{A}(\gamma_2\xi)\frac{d(\mu\circ \gamma^{-1})}{d\mu}(\xi)d\mu(\xi)\\
=&\ \int|\gamma_2'(\xi)|^\delta 1_{A}(\gamma_2\xi)d\mu(\xi).
\end{align*}
 So $\mu_Y$ is defined properly. Since $\mu$ is non-atomic (see Proposition \ref{prop:non_atomic}),  $\mu_Y$ is non-atomic. Moreover, by definition    different choices of  $\mu$  lead to  different measures $\mu_Y$ (consider $\gamma\in \Gamma$).     

We prove that $\mu_Y$ satisfies  Eq. (\ref{eq:mu_Y_gamma_prop}). Given $A\subseteq \partial\mathbb{D}$ Borel and $\gamma_1, \gamma_2\in \Gamma_0$,
\begin{align*}
(\mu_Y\circ \gamma_1)(A\times \{\gamma_2\Gamma     \})=&(\mu_Y )(\gamma_1A\times\{\gamma_1 \gamma_2 \Gamma\})\\
=& \int |(\gamma_1 \gamma_2 )
'(\xi)|^\delta1_{\gamma_1 A}(\gamma_1 \gamma_2 \xi)d\mu(\xi) \\=&
\int |(\gamma_1 \gamma_2 )
'(\xi)|^\delta1_{ A}(\gamma_2 \xi)d\mu(\xi)\\
=& \int |(\gamma_1 )
'(\gamma_2 \xi)|^\delta|(\gamma_2 )
'( \xi)|^\delta1_A(\gamma_2\xi)d\mu(\xi).
\end{align*}  
By the definition of $\mu_Y$,
\begin{align*}
\mu_Y(|\gamma_1'|^\delta 1_{A\times\{\gamma_2\Gamma\}})=\int |(\gamma_1 )
'(\gamma_2 \xi)|^\delta|(\gamma_2 )
'( \xi)|^\delta1_A(\gamma_2\xi)d\mu(\xi)
\end{align*}
and therefore $\frac{d(\mu_Y \circ \gamma_1)}{d\mu_Y}  = |\gamma_1'|^\delta$ for all $\gamma_1\in \Gamma_0$. 

Lastly, we show that this mapping is onto. Given a non-atomic Radon measure $\mu_Y$ which satisfies Eq. (\ref{eq:mu_Y_gamma_prop}),  let 
$\mu(\cdot) := \mu_Y(\cdot, \{\Gamma\})$. Clearly  $\mu$ is non-atomic and  $\mu_Y$ is the resulting measure of the mapping in Eq. (\ref{eq:mu_Y_def}). Moreover, for every $\gamma\in \Gamma$, 
\begin{align*}(\mu\circ \gamma)(A) =&\mu_Y(\gamma A,\{\Gamma\})\\
=&\mu_Y(\gamma A, \{\gamma\Gamma\})\qquad (\because \; \gamma\Gamma=\Gamma \text{ in }G=\Gamma_0/\Gamma)\\
=& (\mu_Y \circ \gamma)(A, \{\Gamma\})\\
= &\mu_Y(|\gamma'|^\delta 1_{A\times\{\Gamma\}})\\
=&\ \int |\gamma'(\xi)|^\delta 1_{A}(\xi)\mu(\xi)
\end{align*} 
and $\mu$ is indeed a  $(\Gamma, \delta)$-conformal measure. 
\hfill\ensuremath{\square}

\section{Conformal measures of hyperbolic covers }
\label{sec:conf_of_hyperbolic_covers}
We now turn out attention to study the conformal measures of a hyperbolic cover. Recall that  $G$ is called a \textit{hyperbolic group} if some (or every, see \cite{ghys_1990}) Cayley graph of $G$ is a hyperbolic graph. We denote by $\partial G$ the hyperbolic boundary of $G$, see definition in Section \ref{sec:hyperbolic_boundary}. We say that a regular cover $\mathbb{D}/\Gamma$  of $\mathbb{D} / \Gamma_0$ is a \textit{hyperbolic cover} if the group of deck transformations $G = \Gamma_0 / \Gamma$ is a hyperbolic group. 

 Our main goal is to prove the following theorem, which describes  the  extremal conformal measures of $\Gamma$ in terms of  $\partial G$:  
\begin{thm}
\label{thm:conv_along_paths}
Let $\Gamma_0$ be a cocompact Fuchsian group, let $ \Gamma\vartriangleleft  \Gamma_{0}$ and let $\delta > \delta(\Gamma)$. Assume that  $G=\Gamma_0 / \Gamma$ is a hyperbolic group.  Then, for every  $\mu \in \confG$,  for $\mu$-a.e. $\xi\in \partial \mathbb{D}$ with Bowen-Series coding $(a_n)$, the  sequence $$e_{a_{n}}^{-1}\dots e_{a_{0}}^{-1} \Gamma$$ converges to a point in $\partial G$. 
If $\mu\in \extconf$, then there exists $\eta\in \partial G$  s.t. the sequence almost-surely converges to $\eta$. Conversely, for every $\eta \in \partial G$, there exists a unique $\mu \in \extconf$ with $\eta$ its almost-surely limiting point of the sequence.
\end{thm}
\begin{remark}\em 
Motivated by the  recent announcement of Bispo and Stadlbaur, we conjecture that the theorem should hold  at the critical value $\delta = \delta(\Gamma)$ as well, if $\mathbf{p}(\Gamma, \delta(\Gamma))<\infty$.
\end{remark}

To prove the theorem, we introduce hyperbolic graph structures on $G$ and $S_X = S_{0}\times G$. Let
$$E_G =\{(g_1, g_2) \in G \times G: g_1=g_2 \text{ or }\exists a\in S_{0} \text{ s.t. } e_{a}^{-1}g_1 = g_2 \} $$
and
let$$E_X =\{((a,g), (b,h)) \in S_X \times S_X: (g,h) \in E_G\}. $$
 Since $\{e_a\}_{a\in S_{0}}$ is a symmetric set  which generates $\Gamma_0$, the set $\{e_a\Gamma\}_{a\in S_{0}}$ generates $G$ and $(G, E_G)$ is an undirected Cayley graph of $G$.
Since $(G, E_G)$ is undirected, $(S_X, E_X)$ is undirected as well. Let $\pi_{S_X}:S_X\rightarrow G$ be the natural projection, $\pi_{S_X}(\xi, g) = g$. Observe that $(S_X, E_X)$ is not the canonical graph associated to the transition matrix of the TMS $X$. In fact, it is  larger and has more edges. 
\begin{defn}\em
Two metric spaces $(X_1, d_1)$ and $(X_2, d_2)$ are called \textit{quasi-isometric} if there exist $g:X_1\rightarrow X_2$, $A\geq 1, B\geq 0$ and $C\geq 0$ s.t.
\begin{enumerate}
\item 
For every $x,y\in X_1$, 
$$\frac{1}{A}d_1(x,y)-B\leq d_2\bigl(g(x), g(y)\bigr)\leq Ad_1(x,y) + B.  $$
\item
For every $y\in X_2$ there exists $x\in X_1$ s.t.
$$d_2\bigl(y, g(x)\bigr)\leq C. $$
\end{enumerate}
We call such a function $g$ a \textit{quasi-isometry}, see \cite{harpe_2000}.
\end{defn}
\begin{prop}
\label{prop:quasi_iso}
The graphs $(G, E_G)$ and $(S_X, E_X)$ are quasi-isometric w.r.t. their natural graph metrics. 
\end{prop}
\begin{proof}
We show that the natural projection $\pi_{S_X}:S_X\rightarrow G$ is a quasi-isometry.
By definition,   $\bigl((a,g), (b,h)\bigr)\in E_X$ iff $(g,h)\in E_G$. Therefore
$$d_{E_G}(g,h) \leq d_{E_X}\bigl((a,g), (b,h)\bigr), \quad \forall (a,g), (b,h)\in S_X. $$
Let $(a_1, g_1), (a_2, g_2)\in S_X$. If $g_1=g_2$ then either $a_1=a_2$ and $d_{E_X}\bigl( (a_1, g_1), (a_2, g_2)\bigr)=0 $ or $a_1\neq a_2$ and $d_{E_X}\bigl(
 (a_1, g_1), (a_2, g_2)\bigr)= 1$. If $d_{E_G}(g_1,g_2)=n>0$, choose $b_1,\dots, b_n\in S_{0}$ s.t. $g_2 = e_{b_n}^{-1}\dots e_{b_1}^{-1} g_1$. By definition, for every $1< i\leq n$, $$\bigl((a_2, e_{b_{i-1}}^{-1}\dots e_{b_1}^{-1} g_1), (a_2, e_{b_i}^{-1}\dots e_{b_1}^{-1} g_1)\bigr)\in E_X  $$
and 
$$\bigl((a_1, g_1), (a_2, a_1^{-1}g_1)\bigr)\in E_X. $$
Therefore $d_{E_X}\bigl(
 (a_1, g_1), (a_2, g_2)\bigr) \leq n$. We conclude that

\begin{equation}
\label{eq:quasi_iso}
d_{E_G}(g_1,g_2)\leq d_{E_X}\bigl( (a_1, g_1), (a_2, g_2)\bigr) \leq d_{E_G}(g_1,g_2) +1\end{equation}
and $\pi_{S_X}$ is indeed a quasi-isometry.
\end{proof}

\begin{cor}
\label{cor:quasi_iso_hyperbolic}
If $G$ is a hyperbolic group, then $(S_X, E_X)$ is a hyperbolic graph.
\end{cor}
\begin{proof}
This follows directly from Proposition \ref{prop:quasi_iso} since hyperbolicity is preserved under quasi-isometries; see \cite{ghys_1990}.
\end{proof}

We denote by $\partial G$ and $\partial S_X$  the hyperbolic boundaries of $(G, E_G)$ and $(S_X, E_X)$ respectively, see definitions is Section \ref{sec:hyperbolic_boundary}. 
\begin{prop}
\label{prop:geo_bound_eq}Assume that $G$ is hyperbolic. Then, the surjection $\pi_{S_X} :S_X \rightarrow G$ extends uniquely to a surjection $\pi_{S_X} :S_X \cup \partial S_X\rightarrow G\cup \partial G$ s.t. $\pi_{S_X}(\partial S_X) =\partial G$ and  
$$(a_n, g_n)\rightarrow \xi \in \partial S_X \Longleftrightarrow g_n\rightarrow \pi_{S_X}(\xi)\in \partial G. $$
In particular, $\pi_{S_X} :\partial S_X \rightarrow \partial G$ is a bijection.
\end{prop}
\begin{proof}
Fix $o_G\in G, o_{S_X}\in S_X$ with $\pi_{S_X}(o_{S_X}) = o_G$.    Recall the definition of $\wedge$ from Section \ref{sec:hyperbolic_boundary}. By Eq. (\ref{eq:quasi_iso}), for every $(a,g), (b,h)\in S_X$, 
\begin{align*}&2|g \wedge h|_{o_{G}}  =   d_{E_G}(o_{G}, g) + d_{E_G}(o_{G}, h) - d_{E_G}(g, h) \\
&\geq  d_{E_X}\bigl(o_{S_X}, (a, g)\bigr) + d_{E_X}\bigl(o_{S_X}, (b, h)\bigr) -d_{E_X}\bigl((a,g), (b, h)\bigr)-2\\
&=2| (a, g) \wedge (b, h)|_{o_{S_{X}}} - 2
\end{align*}
and
\begin{align*}&2|g \wedge h|_{o_{G}}  =   d_{E_G}(o_{G}, g) + d_{E_G}(o_{G}, h) - d_{E_G}(g, h) \\
&\leq  d_{E_X}\bigl(o_{S_X}, (a, g)\bigr) + d_{E_X}\bigl(o_{S_X}, (b, h)\bigr) -d_{E_X}\bigl((a,g), (b, h)\bigr)+1\\
&=2| (a, g) \wedge (b, g)|_{o_{S_{X}}} +1.
\end{align*}
Therefore, 
\begin{equation}
\label{eq:pi_2_wedge}
 |(a, g) \wedge (b, h)|_{o_{S_{X}}}- 1\leq|g\wedge h|_{o_{G}}\leq |(a, g) \wedge (b, h)|_{o_{S_{X}}}+ 1
\end{equation}
and for every $\{a_n\}\subseteq S_{0}$ and $\{g_n\} \subseteq G$, 
\begin{equation*}
 \lim_{m,n\rightarrow\infty}|g_n\wedge g_m|_{o_{G}}= \infty \text{ iff } \lim_{m,n\rightarrow\infty}|(a_n, g_n) \wedge (a_m, g_m)|_{o_{S_{X}}}\rightarrow \infty.
\end{equation*}
In particular, $(a_n,g_n)$ converges to a point in $\partial S_X$ iff $g_n=\pi_{S_X}(a_n,g_n)$ converges to a point in $\partial G$. 

For $\eta\in \partial S_{X}$, set $\pi_{S_X}(\eta) =\lim_{n\rightarrow\infty} g_n$ where $(a_n, g_n)\rightarrow \eta$. Clearly $\pi_{S_X}(\partial S_X)=\partial G$. If $(a_n, g_n)$ and $(b_n, h_n)$ both converges to $\eta\in \partial S_X$, then by Eq. (\ref{eq:pi_2_wedge}) we have that $|g_n\wedge h_n|_{o_G}\rightarrow \infty$ meaning  $h_n$ and $g_n$ both converges to the same limit in $\partial G$. This implies that  $\pi_{S_X}(\xi)$ is well-defined. 

Assume that $\pi_{S_X}(\eta_1) = \pi_{S_X}(\eta_2)$ and let $(a_n, g_n)\rightarrow \eta_1, (b_n, h_n)\rightarrow \eta_2$. Then, $|g_n\wedge h_n|_{o_G}\rightarrow\infty$. By Eq. (\ref{eq:pi_2_wedge}) 
$$|(a_n, g_n)\wedge (b_n, h_n)|_{o_{S_X}}\rightarrow\infty $$
meaning $\eta_1= \eta_2$ and $\pi_{S_X}$ is indeed $1-1$ on the boundary.
\end{proof}
Recall the definition of uniformly irreducibility from Section \ref{sec:Disct_green}.
\begin{prop}
\label{prop:uniformly_irreducible}
$\phi^{X,\delta}$ is uniformly irreducible w.r.t. $(S_{X},E_X)$.
\end{prop} 
\begin{proof}
 Since $\Sigma$ is compact, $\phi^{X,\delta}$ is bounded. If $[(a,g),(b,h)]\neq \varnothing$ then $h=e_a^{-1}g$. In particular $(g,h)\in E_G$ and thus $\bigl((a,g),(b,h)\bigr)\in E_X$. For every $a,b\in S_0$  and $\gamma\in \Gamma_0$, let $n_{a,b,\gamma}$ be  an integer s.t. there is an admissible path from $(a, \Gamma)$ to $(b, \gamma^{}\Gamma)$ in  $X$, namely    
$$\bigl(L_{\phi^{X,\delta}}^{n_{a,b,\gamma}}(1_{[(a,\Gamma)]})\bigr)(x_{(b,\gamma\Gamma)})>0 \quad  $$
where  $x_{(b,\gamma\Gamma)}\in T[(b,\gamma\Gamma)]$. Such a path exists by Proposition \ref{prop:group_ext_transitive}. 

Let $\bigl((a,g),(b,h)\bigr)\in E_X$. Then, either $g=h$ and $$L_{\phi^{X,\delta}}^{n_{a,b,id_{\Gamma_0}}}(1_{[(a,g)]})\bigr)(x_{(b,h)})>0$$ or $h = e_c^{-1} g$, for some $c\in S_0$, and
$$L_{\phi^{X,\delta}}^{n_{a,b,e_c^{-1}}}(1_{[(a,g)]})\bigr)(x_{(b,h)})>0.$$
Thus, with $$K =\max_{a,b\in S_0} \max_{\gamma \in \{e_c\}_{c\in S_0}\cup \{id_{\Gamma_0}\}}n_{a,b,\gamma}$$ we have that $\phi^{X,\delta}$ is uniformly irreducible w.r.t. $(S_X, E_X)$.
\end{proof}

 We are now ready to prove the main result of this section. 
\paragraph{\bf{Proof of Theorem \ref{thm:conv_along_paths}}.}
 By the assumption of the theorem, $\delta > \delta(\Gamma)$. So by Proposition \ref{prop:pressure_delta} we have that $P_G(\phi^{X,\delta})<0$.  By Proposition \ref{prop:uniformly_irreducible} and Corollary \ref{cor:quasi_iso_hyperbolic}, $\phi^{X,\delta}$ is uniformly irreducible w.r.t. the (larger) hyperbolic graph $(S_X, E_X)$. Thus the conditions of Theorem \ref{thm:hyperbolic_main_thm} holds. 

Let $\mu\in \confG$. By Corollary \ref{cor:unique_integral_rep}, we can assume w.l.o.g. that $\mu \in \extconf$.  Let $\mu_X$ be the corresponding eigenmeasure on $X$ from Theorem \ref{prop:gamma_conf_bijection}. Since $\mu$ is extremal and the transformation from $\mu$ to $\mu_X$ is linear, $\mu_X$ is extremal as well. 
By Theorem \ref{thm:hyperbolic_main_thm},  there exists $\eta' \in \partial S_X$ s.t. for $\mu_X$-a.e. $x=(\sigma,  \Gamma)\in X$, $T^nx \rightarrow \eta'$. Let $\eta = \pi_{S_X}(\eta') \in \partial G$. By Proposition \ref{prop:geo_bound_eq}, we have that,   
$$e_{\sigma_n}^{-1}\cdots e_{\sigma_0}^{-1}\Gamma\rightarrow \eta  $$
on the Cayley graph of $G$. Since $\mu(\cdot) =\mu_X\bigl(\pi_{\Sigma}^{-1} (\cdot)\times \{\Gamma\}\bigr )$, $(\sigma, \Gamma)$ is $\mu_X$-typical point iff $\sigma$ is a $\mu$-typical point and thus the first part of the theorem follows.

Now, let $\eta \in \partial G$ and let $\eta' = \pi_{S_X}^{-1}\eta$. By Theorem \ref{thm:hyperbolic_main_thm} there is a unique eigenmeasure $\mu_X$ s.t. for $\mu_X$-a.e. $x\in X$, $T^nx \rightarrow \eta'$. Then, the second part of the theorem follows with $\mu(\cdot) =\mu_X\bigl(\pi_{\Sigma}^{-1} (\cdot)\times \{\Gamma\}\bigr )$.
\hfill\ensuremath{\square}

\section{Convergence of  cutting sequences along geodesics}
\label{sec:conv_geo_flow}
In this section, we study the asymptotic behavior of cutting sequences on hyperbolic covers w.r.t. conformal measures. In particular, for every $\delta > \delta(\Gamma)$ and every $\mu\in \confG$, we show that the cuttings sequence (projected to $G$) $\mu$-a.s. converges to a point in $\partial G$.

We emphasize that the  geodesics on a regular cover do not always escape to infinity. In fact, by the Hopf-Tsuji-Sullivan Theorem the geodesic flow is conservative w.r.t. the Liouville  measure iff the Poincar\'e series diverges for $\delta = 1$, see \cite{aaronson_1997}. An example is a  $\mathbb{Z}^d$-cover: the Poincar\'e series diverges iff $d\leq 2$ \cite{rees_1981}.

In what follows, recall that $F_0 \subseteq \mathbb{D}$ is a  fundamental domain of $\mathbb{D} / \Gamma_0$. For every $\gamma_1,\gamma_2\in \Gamma_0$, $$int(\gamma_1 F_0) \cap int(\gamma_2 F_0) \neq \varnothing\Longleftrightarrow\ \gamma_1 = \gamma_2$$ and  $$\gamma_1 F_0 \text{ and }\gamma_2 F_0 \text{ share a common edge } \Longleftrightarrow\ \gamma_1 \gamma_2^{-1}\in \{e_a\}_{a\in S_0}.$$Given $\xi^-, \xi^+\in \partial \mathbb{D}$ with $\xi^- \neq \xi^+$, we denote by $\xi^- \wedge \xi^+$ the unique geodesic curve in $\mathbb{D}$ from $\xi^-$ to $\xi^+$. Let 
$$\mathcal{R} = \bigl\{(\xi^- ,\ \xi^+)\in \partial\mathbb{D}^2   : (\xi^- \wedge\ \xi^+)   \cap int(F_{0}) \neq \varnothing\bigr\}. $$  
Observe that $\mathcal{R}$ is symmetric, namely  $(\xi^-, \xi^+)\in \mathcal{R}$ if and only if $(\xi^+, \xi^-)\in \mathcal{R}$. 
Let $(\xi^-, \xi^+)\in \mathcal{R}$ and let $\{F_{i}\}_{i\in \mathbb{Z}}$ be the sequence of copies of $F_0$ that the curve $(\xi^-\wedge  \xi^+)$ intersects. In case $(\xi^-\wedge  \xi^+)$ passes through a vertex of some $F_i$, we perturb the curve around it, see Figure 5 in \cite{series_1986}.   Then, for all $i$ there exists a unique $e_i\in \{e_a\}_{a\in S_0}$ s.t. $ F_i = e_i^{-1} F_{i+1}$.
\begin{defn}\em
The sequence $(\dots, e_{{-1}}, e_{0}, e_{1},\dots)$ is called the \textit{cutting sequence} of $(\xi^-, \xi^+)$.  
\end{defn} 

\begin{defn}\em
For $(\xi^-, \xi^+) \in \mathcal{R}$ with a cutting sequence $(\dots, e_{{-1}}, e_{0}, e_{1},\dots)$, we write
$$\eta_n^+(\xi^-, \xi^+) =e_{n}^{-1}\cdots e_0^{-1} \Gamma$$ 
and
$$\eta_{n}^-(\xi^-, \xi^+) =(e_{-n-1})^{-1}\cdots (e_{-1})^{-1} \Gamma.$$
\end{defn}
The following theorem describes the the limiting behaviour of $\eta_n^+$ and $\eta_n^-$ w.r.t. a conformal measure.
\begin{thm}
\label{thm:conv_cutting_seq}
Let $\Gamma_0$ be a cocompact Fuchsian group, let $ \Gamma\vartriangleleft  \Gamma_{0}$, let $\delta > \delta(\Gamma)$ and let $\mu\in \confG$. Assume that $G=\Gamma/\Gamma_0$  is a hyperbolic group.   Then, 
\begin{enumerate}
\item 
For $ \mu$-a.e. $\xi^+\in\partial \mathbb{D}$, for every $\xi^-\in\mathbb{D}$ s.t. $(\xi^-, \xi^+)\in \mathcal{R}$, the sequence $\eta_{n}^+(\xi^-, \xi^+)$ converges to a point in $\partial G$.  
\item 
For $ \mu$-a.e. $\xi^-\in\partial \mathbb{D}$, for every $\xi^+\in\mathbb{D}$ s.t. $(\xi^-, \xi^+)\in \mathcal{R}$, the sequence $\eta_{n}^-(\xi^-, \xi^+)$ converges   to a point in  $\partial G$.  

\end{enumerate}
If $\mu\in \extconf$  then there exists a point $\eta\in \partial G$ s.t. the sequences  almost-surely converges to $\eta$. Conversely, for every $\eta \in \partial G$, there exists a unique $\mu \in \extconf$ with $\eta$ its almost-surely limiting point of the sequences. 
\end{thm}
\begin{remark}\em
The limiting point $\eta$ is the same limiting point from Theorem \ref{thm:conv_along_paths}. Again, motivated by the recent announcement of Bispo and Stadlbaur, we conjecture that the theorem  should hold  at the critical value $\delta = \delta(\Gamma)$, if $\mathbf{p}(\Gamma, \delta(\Gamma))<\infty$.
\end{remark}

To prove the theorem, we  exploit the connection between  boundary expansions and  cutting sequences, presented by Series \cite{series_1986}.
To do so, we briefly introduce the two-sided Bowen-Series coding. Denote by $\Sigma^+=\Sigma$  the \textit{positive one-sided shift} (see Eq. (\ref{eq:pos_sidede_shift})),
let $$\Sigma^-=\bigl\{(\dots, \sigma_{-2}, \sigma_{-1}, \sigma_0): \forall i<0, [\sigma_i, \sigma_{i+1}]\neq \varnothing \text{ in }\Sigma^+\bigr\}$$ be the \textit{negative one-sided shift} and let $$\Sigma^\pm=\bigl\{(\dots,  \sigma_{-1}, \sigma_0, \sigma_1,\dots): \forall i, [\sigma_i, \sigma_{i+1}]\neq \varnothing \text{ in }\Sigma^+\bigr\} $$be the  \textit{two-sided shift}. We write $T_\Sigma$ for the left-shift action  both on $\Sigma^+$ and $\Sigma^{\pm}$; the meaning should be clear from the context. Recall that
 $\pi_\Sigma:\Sigma^{+}\rightarrow\partial\mathbb{D}$ is the canonical projection where $$\pi_\Sigma(\sigma_0, \sigma_1, \dots) \in \bigcap_{n\geq 0}\overline{f_{\Gamma_0}^{-n}I_{\sigma_n}}$$ We write $\pi^+:\Sigma ^{\pm}\rightarrow\partial\mathbb{D}$, $$\pi^+(\sigma)= \pi_\Sigma(\sigma_0, \sigma_1, \dots).$$
For every $a\in S_0$, let $\overline a \in S_0$ s.t. $e_{\overline a} = e_a^{-1}$ and $[a,b]\neq \varnothing$ iff $[\bar b, \overline a]\neq \varnothing$ in $\Sigma^+$.  We define$$\pi^-(\sigma) = \pi_\Sigma(\overline {\sigma_{-1}}, \overline {\sigma_{-2}},\dots) $$
and
$$\pi(\sigma) = (\pi^-(\sigma), \pi^+(\sigma)). $$
Then, $$\pi(T_{\Sigma}\sigma) = \bigl(\pi^-(T_{\Sigma}\sigma),\pi^+(T_{\Sigma}\sigma)\bigr) =\bigl (e_{\sigma_0}^{-1}\pi^-(\sigma), e_{\sigma_0}^{-1}\pi^+(\sigma)\bigr).$$

Let
$$\mathcal{A} =\bigl \{(\xi^- ,\ \xi^+)  : \exists \sigma \in \Sigma^{\pm} \text{ s.t. } (\xi^-, \xi^+)=\pi(\sigma)\bigr\}. $$
The Bowen-Series map $f_{\Gamma_0}$ acts on $\mathcal{A}$ similarly to the left-shift action,
$$f_{\Gamma_0}(\xi^-,\ \xi^+) =(\pi \circ T_\Sigma)(\sigma)=\bigl(e_{\sigma_0}^{-1}\pi^-(\sigma),e_{\sigma_0}^{-1}\pi^+(\sigma)\bigr) = e_{\sigma_0}^{-1}(\xi^- ,\ \xi^+) . $$
  Here and throughout $e_{\sigma_0}^{-1}(\xi^- ,\ \xi^+):=(e_{\sigma_0}^{-1}(\xi^- ),\ e_{\sigma_0}^{-1}(\xi^+))$. The value of the  \textit{first-return map} $g_{\Gamma_0}:\mathcal{R}_{}\rightarrow\mathcal{R}$ on a pair $(\xi^- ,\ \xi^+) \in \mathcal{R}$ with  cutting sequence $(\dots,e_{{-1}},e_{0},e_{1},\dots)$ is $$g_{\Gamma_0} (\xi^- ,\ \xi^+)= e_{0}^{-1}(\xi^- ,\ \xi^+).$$
Notice that $g_{\Gamma_0} (\xi^- ,\ \xi^+) \in \mathcal{R}$  and the cutting sequence of $g_{\Gamma_0} (\xi^- ,\ \xi^+) $ is the cutting sequence of $(\xi^- ,\ \xi^+)$ shifted by one position to the left. In particular,  
$$g^{n}_{\Gamma_0} (\xi^- ,\ \xi^+)= (e_{{n-1}}^{-1}\cdots e_{0}^{-1})(\xi^- ,\ \xi^+). $$

\begin{thm}[Series \cite{series_1986}]
\label{thm:series_cutting_coding}
There exists a bijection $\varphi :\mathcal{A}\rightarrow\mathcal{R}$ s.t.
$$ (\varphi \circ f_{\Gamma_0})(\xi^- ,\ \xi^+) =( g_{\Gamma_0} \circ \varphi )(\xi^- ,\ \xi^+). $$
\end{thm}

For a given $\gamma_0\in \Gamma_0$, we denote by $|\gamma_0|$ the \textit{word length} of $\gamma_0$:  the minimal integer s.t. there exist admissible $a_1,\dots, a_{|\gamma_0|}\in  S_{0}$ with $\gamma_0 = e_{a_1}\cdots e_{a_{|\gamma_0|}}$. We let $|id_{\Gamma_0}|=0$. The following lemma  shows that the action of the bijection $\phi$ is uniformly bounded.

\begin{figure}[ht!]
\centering
\includegraphics[width=2.0in,height=2.0in]{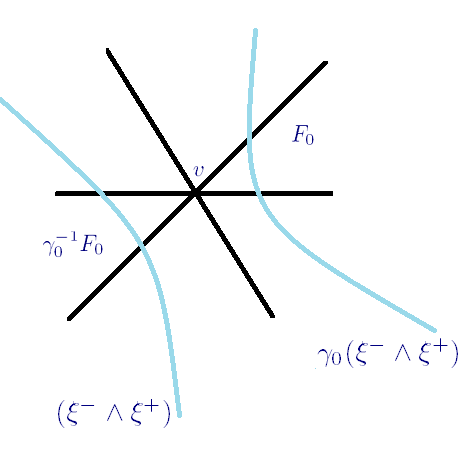}
\caption{An illustration of  $F_0$ and $\gamma_0^{-1}F_0$ sharing a common vertex $v$. In this figure, $n(v) = 6$.}
\label{fig:sharing_vertex}
\end{figure}

\begin{lemma}
\label{lemma:series_bounded}
There exists $n_{\Gamma_0}>0$ s.t. for every $(\xi^-, \xi^+)\in \mathcal{A}$, there is $\gamma_0\in \Gamma_0$ with  $|\gamma_0| \leq n_{\Gamma_0}$ and  $\varphi(\xi^- ,\ \xi^+) = \gamma_0 (\xi^- ,\ \xi^+)$.
\end{lemma}
\begin{proof}
 Let $(\xi^-, \xi^+)\in \mathcal{A}$. To transform the curve $(\xi^-\wedge \xi^+)$ to a curve that intersects $int(F_0)$, Series paired $F_0$ with a different copy $\gamma_0^{-1} F_0$ where $(\xi^-\wedge \xi^+) \in int(\gamma_0^{-1} F_0)$. Then, $ \gamma_0(\xi^-\wedge \xi^+)\in \mathcal{R}$. When $(\xi^-, \xi^+)\in \mathcal{R}$, $\gamma_0 = id_{\Gamma_0}$ and $|\gamma_0|=0$.  

Assume that $(\xi^-, \xi^+)\not\in \mathcal{R}$. By Proposition 3.2 in \cite{series_1986}, there is no additional copy of $F_0$ between $F_0$ and $\gamma_0^{-1} F_0$. In particular, there are two possible scenarios: either $F_0$ and $\gamma_0^{-1} F_0$ share a common edge or either they share a single vertex. 

If they share a common edge then $|\gamma_0|=1$. Assume that they share exactly one vertex $v$. See  Figure \ref{fig:sharing_vertex}. Let $n(v)$ be the degree of $v$ and let $\gamma_v^1, \dots, \gamma_{v}^{n(v)}$ be the transformations between $F_0$ to the adjoint copies of $F_0$ that share the vertex $v$ with $F_0$. Then,   $|\gamma_0|\leq \max_{i\leq n(v)}|\gamma_v^i|$. Since $\Gamma_0$ is co-compact, $\partial F_0$ has finite number of  vertices, all with finite degrees, and   the lemma follows with  $n_{\Gamma_0}= \max_{v\in \partial F_0}\max_{i\leq n(v)}|\gamma_v^i|$.  
\end{proof}

 Let $$\mathcal{D} := \{\xi \in \partial\mathbb{D}\setminus Fix(\Gamma_0) : |\pi_{\Sigma}^{-1}(\xi)|=1\}.$$ Notice that $\mathcal{D}$ is a $\Gamma_0$-invariant set and that $\partial \mathbb{D} \setminus \mathcal{D}$ is countable and thus a null set w.r.t. any conformal measure.  Given $\xi \in \mathcal{D}$ with a (one-sided) boundary expansion $\pi_{\Sigma}^{-1}(\xi)=(\sigma_0,\sigma_1 \dots)$ and $\gamma_0\in \Gamma_0$, we write
$$\tau_n(\xi, \gamma_0) :=e_{\sigma_{n}}^{-1}\cdots e_{\sigma_{0}}^{-1}\gamma_{0}\Gamma.$$
To prove Theorem \ref{thm:conv_cutting_seq}, we first introduce and prove   two auxiliary lemmas.

\begin{lemma}
\label{lemma:dist_eta_tau}
Let $\xi^+\in \mathcal{D}$ and let $\gamma_0\in \Gamma_0$. Then, for all $\xi^-\in \mathcal{D}$  with $(\xi^-, \xi^+)\in \mathcal{R}$ and $\varphi^{-1}(\xi^-, \xi^+) = \gamma_0(\xi^-, \xi^+)$, 
$$d_{E_G}\bigl(\eta_{n}^+(\xi^-, \xi^+),\tau_n(\gamma_0\xi^{+}, \gamma_0)\bigr)\leq n_{\Gamma_0}  $$ 
 where $n_{\Gamma_0}$ is the bound from Lemma \ref{lemma:series_bounded}.
\end{lemma}
\begin{proof}Let $\xi^-\in \mathcal{D}$ s.t. $(\xi^-, \xi^+)\in \mathcal{R}$ and $\varphi^{-1}(\xi^-, \xi^+) = \gamma_0(\xi^-, \xi^+)$. Let $\sigma\in \Sigma^+$ s.t. $\pi_{\Sigma}(\sigma) = \gamma_0\xi^{+}$ and let $(\dots, e_{-1}, e_0, e_1, \dots)$ be the cutting-sequence of $(\xi^-, \xi^+)$.  For every $n\geq 1$, let $\gamma_n \in \Gamma_0$ s.t. $\varphi^{-1}( g^{n}_{\Gamma_0}(\xi^-, \xi^+)) =\gamma_n (g^{n}_{\Gamma_0}(\xi^-, \xi^+))$.  By Theorem \ref{thm:series_cutting_coding}, 
$$ (\varphi ^{}\circ f^{n}_{\Gamma_0}\circ \varphi^{-1})(\xi^- ,\xi^+) = g^{n}_{\Gamma_0} (\xi^- ,\ \xi^+)$$ 
meaning 
$$(\gamma^{-1}_n e_{\sigma_{n-1}}^{-1}\cdots e_{\sigma_{0}}^{-1}\gamma_0^{})(\xi^- ,\ \xi^+)=  (e_{{n-1}}^{-1}\cdots ^{}e_{0}^{-1})(\xi^- ,\ \xi^+).$$ 
Since $\xi^-, \xi^+\not\in Fix(\Gamma_0)$, 
$$\gamma^{-1}_n e_{\sigma_{n-1}}^{-1}\cdots e_{\sigma_{0}}^{-1}\gamma_0^{} =e_{{n-1}}^{-1}\cdots ^{}e_{0}^{-1}.  $$
By Lemma \ref{lemma:series_bounded}, $|\gamma_n|\leq n_{\Gamma_0}$ and therefore   \begin{align*}&d_{E_G}\bigl(\eta_{n}^+(\xi^-, \xi^+),\tau_n(\gamma_0\xi^{+}, \gamma_0)\bigr)\\
&=d_{E_G}(\gamma^{-1}_n e_{\sigma_{n-1}}^{-1}\cdots e_{\sigma_{0}}^{-1}\gamma_0\Gamma, ^{} e_{\sigma_{n-1}}^{-1}\cdots e_{\sigma_{0}}^{-1}\gamma_0\Gamma)\leq n_{\Gamma_0}.\end{align*}
\end{proof}
\begin{lemma}
\label{lemma:eta_plusminus}
Let $(\xi^-, \xi^+)\in \mathcal{R}$. Then, 
$$\eta_n^+(\xi^+, \xi^-) = \eta_n^-(\xi^-, \xi^+).  $$
\end{lemma}
\begin{proof}
Observe that if  the cutting sequence of $(\xi^-, \xi^+)$ is $(e_i)$ then the cutting sequence of $(\xi^+, \xi^-)$ is $(f_i)$ with $f_i=e_{-i-1}$.  This implies that $$\eta_{n}^+(\xi^+, \xi^-) ={(f_{n})}^{-1}\cdots {(f_{0})}^{-1} \Gamma_{} ={(e_{-n-1})}^{-1}\cdots {(e_{-1})}^{-1} \Gamma = \eta_{n}^-(\xi^-, \xi^+).$$
\end{proof}
\paragraph{\bf{Proof of Theorem \ref{thm:conv_cutting_seq}.}} 
We show that for $\mu$-a.e. $\xi^+\in \partial \mathbb{D}$, for every  $\xi^-\in \mathcal{D}$ with $(\xi^-, \xi^+)\in \mathcal{R}$, $\eta_n^+(\xi^-, \xi^+)$ converges to a point in $\partial G$. By Lemma \ref{lemma:eta_plusminus}, the arguments for $\eta_{n}^-(\xi^-, \xi^+)$ are similar. By Corollary \ref{cor:unique_integral_rep}, we can assume w.l.o.g. that $\mu^+$ is extremal. Since $\mu$ is non atomic (see Proposition \ref{prop:non_atomic}), we can also assume that $\xi^-, \xi^+ \in \mathcal{D}$. 

Let $\gamma_0\in \Gamma_0$, let $\eta \in \partial G$ and let
$$A_{\gamma_0} =\left\{  \xi  ^{+}\in \mathcal{D}: \begin{matrix}\exists \xi^-\in \mathcal{D} \text{ s.t. }(\xi^-, \xi^+)\in \mathcal{R}, \;\varphi^{-1}(\xi^-, \xi^+) = \gamma_0(\xi^-, \xi^+),  \\
 \text{and } \lim_{n\rightarrow\infty}\eta_{n}^+(\xi^-, \xi^+)\neq \eta   \\
\end{matrix}\right\}. $$
We write $\lim_{n\rightarrow\infty}\eta_{n}^+(\xi^-, \xi^+)\neq \eta$ whenever the limit does not exist or it exists but differs from $\eta$. We show  that there exists $\eta\in \partial G$ s.t. $\mu(A_{\gamma_0})=0$.

  Given  $\xi^+\in A_{\gamma_0}$  and  $\xi^-\in \mathcal{D}$ s.t. $(\xi^-, \xi^+)\in \mathcal{R}$ and $\varphi^{-1}(\xi^-, \xi^+)=\gamma_0(\xi^-, \xi^+)$, we have   by Lemma \ref{lemma:dist_eta_tau} that 
\begin{equation*}
\label{eq:convergence_condition}
\lim_{n\rightarrow\infty}\eta_{n}^+(\xi^-, \xi^+) = \eta \Longleftrightarrow\  \lim_{n\rightarrow\infty} \tau_n(\gamma_0\xi^{+}, \gamma_0) =\eta. 
\end{equation*}
Observe that the right term does not depend on $\xi^-$ once $\gamma_0$ and $\xi^+$ are given. 
Hence,
$$A_{\gamma_0}=\bigl\{\xi  ^{+}\in \mathcal{D}:   \lim_{n\rightarrow\infty} \tau_n(\gamma_0\xi^{+}, \gamma_0) \neq\eta\bigr\}. $$
Since $\mathcal{D}$ is a $\Gamma_0$-invariant set,
$$\gamma_0 A_{\gamma_0} = \bigl\{\xi\in \mathcal{D}  : \lim_{n\rightarrow\infty} \tau_n(\xi, \gamma_0) \neq\eta\bigr\}.$$
Let $\mu_Y$ the measure on $Y=\partial \mathbb{D}\times G$ from Theorem \ref{prop:gamma_conf_bijection} that satisfies Eq. (\ref{eq:mu_Y_gamma_prop}) and
$$\mu(\cdot) = \mu_Y(\cdot\times \{\Gamma\}). $$
Then, $\mu(A_{\gamma_0}) =\mu_Y(A_{\gamma_0}\times\{\Gamma\})$. By Eq. (\ref{eq:mu_Y_gamma_prop}),   
$$  \mu_Y(A_{\gamma_0}\times\{\Gamma\})=0 \Longleftrightarrow  \mu_Y(\gamma_0A_{\gamma_0}\times\{\gamma_0 \Gamma\})=0.  $$
Let $\mu_X = \mu_Y \circ \tilde \pi$, see Lemma \ref{lemma:mu_X_mu_Y_bijection}. Then, 
$$\mu_Y(\gamma_0A_{\gamma_0}\times\{\gamma_0 \Gamma\})=0 \iff \mu_X\bigl(\pi_\Sigma^{-1}(\gamma_0 A_{\gamma_0} )\times \{\gamma_0 \Gamma\}\bigr)=0. \ $$
By Proposition \ref{prop:geo_bound_eq}, 
$$\pi_\Sigma^{-1}(\gamma_0 A_{\gamma_0} ) = \biggl\{\sigma^+ \in  \pi_\Sigma^{-1}\bigl( \mathcal{D}\bigr  ): \lim_{n\rightarrow\infty}T^{n}(\sigma, \gamma_0)\neq \pi_{S_X}^{-1}(\eta)\biggr\}. $$ 
Then, by Theorem \ref{thm:hyperbolic_main_thm} there exists $\eta\in \partial G$\ s.t. $\pi_\Sigma^{-1}(\gamma_0 A_{\gamma_0} )\times \{\gamma_0 \Gamma\}$ is a $\mu_X$-null set.

Similarly, given $\eta\in \partial G$, let $\mu_X$ s.t. $\pi_\Sigma^{-1}(\gamma_0 A_{\gamma_0} )\times \{\gamma_0 \Gamma\}$ is a $\mu_X$-null set. Such a measure exists by Theorem \ref{thm:hyperbolic_main_thm}. Let $\mu\in \extconf$ with $$\mu(\cdot ) = \mu_X (\pi_\Sigma^{-1}(\cdot) \times \{\Gamma\}). $$
Then, similar arguments show that $\mu(A_{\gamma_0})=0$. 
$$\mu $$       
\hfill\ensuremath{\square}

%% file: ancona.tex
\section{The Martin boundary of a Markov shift over a hyperbolic graph}
\subsection{Discretized Green's function and related inequalities}
\label{sec:Disct_green}
 To prove Theorem \ref{thm:hyperbolic_main_thm}, we introduce approximated versions of the Green's function and the Martin kernel to the discrete set of states $S$ rather than the non-discrete set of infinite paths $X$. For these discretized functions, we present several combinatorial inequalities, inspired by their probabilistic analogues.

 The following proposition shows that observing the first coordinate alone  suffices to determine whether a sequence of internal points $x^n\in X$ converges to a boundary point $\omega\in \mathcal{M}_m(\lambda)$.

\begin{prop}
\label{prop:FirstCoordConvergence}
Assume that $(X,T)$ is locally compact and topologically transitive and that $\phi$  is $\lambda$-transient potential with summable variations. Let $x^n \in X$ with $x^n \xrightarrow[n\rightarrow\infty]{} \omega\in \mathcal{M}_m(\lambda)$ and let $y^{n}\in X$ s.t. $y^{n}_0 = x^{n}_0$ for all $n>0$. Then, $y^{n} \xrightarrow[n\rightarrow\infty]{}\omega$ as well. 
\end{prop}
\begin{proof}
Assume that $y^{n}\rightarrow \omega'$, otherwise we can take a converging sub-sequence. Since $\phi$ has summable variations, for every $f$ of the form $f = 1_{[a_1,\dots, a_n]}$ and for every $n>0$,
$$ G(f, x^n|\lambda) = C_\phi^{\pm 1}G(f, y^{n}|\lambda)$$
where $C_\phi =\exp \bigl( \sum_{k\geq 2}Var_k(\phi)\bigr)$. In particular, 
$$ |\log K(f, x^n|\lambda) - \log K(f, y^{n}|\lambda)| \leq 2 \log C_\phi.$$
By taking $n\rightarrow\infty$, we obtain that
\begin{equation}
\label{eq:omega_omegatag}
 (C_\phi)^{-2} \mu_{\omega'}(f)\leq \mu_\omega (f)\leq C_\phi ^{2}\mu_{\omega'}(f).
 \end{equation}
Since the collection of indicators of cylinder sets linearly spans a dense subset of $C_c(X)$, the inequality in Eq. (\ref{eq:omega_omegatag}) holds for all $f\in C_c^+(X)$. 
 Since $\mu_\omega$ is minimal and $\mu_\omega([o]) = \mu_{\omega'}([o])=1$ , we have that $\mu_\omega = \mu_\omega'$ and thus $\omega = \omega'$.
\end{proof}
\begin{cor}
\label{cor:conv_seq_states}
Under the assumptions of Proposition \ref{prop:FirstCoordConvergence}, for every $\omega \in \mathcal{M}_m(\lambda)$ there exists a sequence $a_n\in S$ s.t. for every sequence  $x^n\in X$ with $x^n\in [a_n]$, 
$$\lim_{n\rightarrow\infty}K(f, x^n|\lambda)=K(f,\omega|\lambda) , \quad \forall f\in C_c(X). $$
\end{cor}

Recall that for every $a\in S$ we fixed $x_a\in T[a]$ arbitrarily. 
\begin{defn}\em
For $a,b\in S$ and $\lambda>0$, let
$$ G(a,b|\lambda): = G(1_{[a]}, bx_b|\lambda) = \sum_{n=0}^\infty\sum_{\substack{(a_{0},\dots,a_n)\in \mathcal{W}^{n+1}\\ a_{0}=a, a_n=b}}\lambda^{-n}e^{\phi_n(a_{0},\dots,a_nx_{b})} $$
and let
$$ F(a,b|\lambda) :=  \sum_{n=0}^\infty\sum_{\substack{(a_{0},\dots,a_n)\in \mathcal{W}^{n+1}\\ a_{0}=a, a_n=b\\ \forall i<n,a_i\neq b}}\lambda^{-n}e^{\phi_n(a_{0},\dots,a_nx_{b})} . $$
We let $\phi_0\equiv0$. In particular,  $F(a,a|\lambda)=1$. 

For a subset $A \subseteq S$, we let
$$ L^A(a,b|\lambda) :=\sum_{n=0}^\infty\sum_{\substack{(a_{0},\dots,a_n)\in \mathcal{W}^{n+1}\\ a_{0}=a, a_n=b\\ a_0 \in A,\;\forall i>0: a_i\not\in A}}\lambda^{-n}e^{\phi_n(a_{0},\dots,a_nx_{b})} $$
and let
$$ F^A(a,b|\lambda) :=\sum_{n=0}^\infty\sum_{\substack{(a_{0},\dots,a_n)\in \mathcal{W}^{n+1}\\ a_{0}=a, a_n=b\\a_n\in A,\; \forall i<n:a_i\not\in A\\}}\lambda^{-n}e^{\phi_n(a_{0},\dots,a_nx_{b})}. $$
Observe that if $a\not\in A$ then $L^A(a,b|\lambda)=0$ and if $b\not\in A$ then $F^A(a,b|\lambda)=0$. We  write $L^{\{a\}}(a,b|\lambda) = L(a,b|\lambda)$.  Let$$K(a,b|\lambda) :=K(1_{[a]}, bx_b|\lambda) =  \frac{G(a,b|\lambda)}{G(o,b|\lambda)}$$ and given $f\in C_c(X)$, let
$$K(f,a|\lambda) := K(f, ax_b|\lambda). $$
For $\lambda=1$, we simply write $G(a,b), F(a,b), L(a,b)$ and $K(a,b)$. 
\end{defn}

The following propositions present several useful inequalities involving the functions $F,G,L$ and $K$. Their proofs are elementary and  included here for completeness.  Several of these inequalities have been adapted from the probabilistic settings; see \cite{woess_2000} for more details.

\begin{prop}
\label{prop:green_main_prop}
Assume that $(X,T)$ is locally compact and transitive, that $\phi$ has summable variations and that $P_G(\phi)<\infty$. Then, there exist a  constant $C>1$ s.t. for every $\lambda> \rho(\phi)$, 
\begin{enumerate}[label={(\arabic*)}]
\item 
\label{prop:green_main_prop_FG}
For every $a,b\in S$, 
$$ G(a,b|\lambda) = C^{\pm 1}F(a,b|\lambda)G(b,b|\lambda).$$
\item 
\label{prop:green_main_prop_F}
For every $a,b,c\in S$, 
$$F(a,c|\lambda)F(c,b|\lambda) \leq CF(a,b|\lambda).$$
\item
\label{prop:green_main_prop_GL}
For every $a,b\in S$ and  every set  $A\subseteq S$ s.t. every path from  $a$ to $b$ must pass through $A$,
$$ G(a,b|\lambda) = C^{\pm 1}\sum_{e\in A}G(a,e|t)L^A(e,b|\lambda).$$
\item
\label{prop:green_main_prop_GLFG}
For every $a,b\in S$ and every set $A\subseteq S$,
$$\sum_{e\in A}G(a,e|\lambda)L^A(e,b|\lambda) = C^{\pm1}\sum_{e\in A}F^A(a,e|\lambda)G(e,b|\lambda).$$

\item
\label{prop:green_main_prop_FEG}
For every $a,b\in S$ and $A\subseteq S$, 
$$\sum_{e\in A}F^A(a,e|\lambda)G(e,b|\lambda)\leq C G(a,b|\lambda).$$
\item 
\label{prop:green_main_prop_Gt}
For every $a,b\in S$ and every $\lambda_1,\lambda_2$ with $\rho(\phi) < \lambda_1 \leq  \lambda_2$, 
$$\frac{G(a,b|\lambda_{1})}{\lambda_1} - \frac{G(a,b|\lambda_2)}{\lambda_2} = C^{\pm 1}\left(\frac{1}{\lambda_1} - \frac{1}{\lambda_2}\right)\sum_{c\in S}G(a,c|\lambda_1)G(c,b|\lambda_2).$$
\item
\label{prop:green_main_prop_K}
For every admissible $a_1,\dots, a_N\in S$ and every $b_1,b_2
\in S$ with $b_i\neq a_j$,
$$ K(1_{[a_1,\dots, a_N]},b_i|\lambda) = C^{\pm 1} t^{-(N-1)}e^{\phi_{N-1}(a_1,\dots, a_{N}x_{a_N})}K(a_{N},b_i|\lambda) $$
and
$$ \frac{K(1_{[a_1,\dots, a_N]},b_1|\lambda)}{K(1_{[a_1,\dots, a_N]},b_2|\lambda)}=C^{\pm 1} \frac{F(a_N, b_1|t)F(o,b_2|\lambda)}{F(a_N, b_2|t)F(o, b_1|\lambda)}.$$
\end{enumerate}
\end{prop}
\begin{proof} 
\begin{enumerate}[label={(\arabic*)}]
\item
\begin{align*}
&G(a,b|\lambda) = \sum_{n=0}^{\infty}\sum_{\substack{(a_0,\dots,a_n)\in \mathcal{W}^{n+1}\\
a_0=a, a_n=b }}\lambda^{-n}e^{\phi_n(a_{0}\dots a_{n}x_b) }\\
=&  \sum_{n=0}^{\infty}\sum_{k=0}^{n}\sum_{\substack{(a_0,\dots,a_n)\in \mathcal{W}^{n+1}\\a_0=a, a_n=b\\ a_k=b, \; \forall i < k:a_i\neq b }}\lambda^{-n}e^{\phi_n(a_{0}\dots a_{n}x_b) }\\
\asymp & \sum_{n=0}^{\infty}\sum_{k=0}^{n}\sum_{\substack{(a_0,\dots,a_n)\in \mathcal{W}^{n+1}\\a_0=a, a_n=b,\; \\
a_k = b, \;\forall i < k:a_i\neq b }}\lambda^{-k}e^{\phi_k(a_{0}\dots a_{k}x_b) }\lambda^{-(n-k)}e^{\phi_{n-k}(a_{k},\dots, a_n x_b)}\\
=& \left(\sum_{k=0}^{\infty}\sum_{\substack{(a_0,\dots,a_k)\in \mathcal{W}^{k+1}\\a_0=a,\;a_k=b,\; \forall i < k:a_i\neq b }}\lambda^{-k}e^{\phi_k(a_{0}\dots a_{k}x_b) }\right)\times \left(\sum_{n=0}^{\infty}\sum_{\substack{(a_0,\dots,a_n)\in \mathcal{W}^{n+1}\\a_0=a_n=b\\  }}\lambda^{-n}e^{\phi_n(a_{0}\dots a_{n}x_b) }\right)\\
 =& F(a,b|\lambda)G(b,b|\lambda).
\end{align*}
Moreover, all estimates are uniform in $a$ and $b$. 
\item
First, observe that
\begin{align*}
& G(a,b|\lambda) \geq \sum_{n=0}^\infty\sum_{\substack{(a_0, \dots, a_n)\in \mathcal{W}^{n+1}\\a_0 = a, a_n = b\\ \exists 0\leq i\leq n: a_i=c}}\lambda^{-n}e^{\phi_n(a_0\dots a_n x_b)}\\
=& \sum_{n=0}^\infty\sum_{k=0}^n\sum_{\substack{(a_0, \dots, a_n)\in \mathcal{W}^{n+1}\\a_0 = a, a_n = b\\  a_k=c, \forall i<k: a_i \neq c}}\lambda^{-n}e^{\phi_n(a_0\dots a_n x_b)}\\
\asymp&\left(\sum_{k=0}^{\infty}\sum_{\substack{(a_0,\dots,a_k)\in \mathcal{W}^{k+1}\\a_0=a,\;a_k=c,\; \forall i < k:a_i\neq c }}\lambda^{-k}e^{\phi_k(a_{0}\dots a_{k}x_c) }\right)\times \left(\sum_{n=0}^{\infty}\sum_{\substack{(a_0,\dots,a_n)\in \mathcal{W}^{n+1}\\a_0=c,a_n=b\\  }}\lambda^{-n}e^{\phi_n(a_{0}\dots a_{n}x_b) }\right)\\ 
 =& F(a,c|\lambda) G(c,b|\lambda).
\end{align*}
Therefore $G(a,b|\lambda) \gg F(a,c|\lambda) G(c,b|\lambda)$. 
By Proposition \ref{prop:green_main_prop}.\ref{prop:green_main_prop_FG},
$$
F(a,b|\lambda) \asymp \frac{G(a,b|\lambda)}{G(b,b|\lambda)} \gg \frac{F(a,c|\lambda)G(c,b|\lambda)}{G(b,b|\lambda)} \asymp F(a,c|\lambda) F(c,b|\lambda). 
$$
Moreover, all estimates are uniform in $a,b$ and $c$.
\item
 Assume that every path from $a$ to $b$ must pass through the set $A\subseteq S$. Then, 
\begin{align*}
&G(a,b|\lambda) = \sum_{n=0}^{\infty}\sum_{k=0}^{n}\sum_{\substack{(a_0,\dots,a_n)\in \mathcal{W}^{n+1}\\a_0=a, a_n=b\\ a_k \in A, \; \forall i>k:a_i\not\in A\\  }}\lambda^{-n}e^{\phi_n(a_{0}\dots a_{n}x_b) }\\
=& \sum_{n=0}^{\infty}\sum_{k=0}^{n}\sum_{e\in A}\sum_{\substack{(a_0,\dots,a_n)\in \mathcal{W}^{n+1}\\a_0=a, a_n=b\\ a_k = e, \; \forall i>k:a_i\not\in A }}\lambda^{-n}e^{\phi_n(a_{0}\dots a_{n}x_b) }\\
\asymp& \sum_{e\in A}\left(\sum_{k=0}^{\infty}\sum_{\substack{(a_0,\dots,a_k)\in \mathcal{W}^{k+1}\\a_0=a, a_k=e }}\lambda^{-k}e^{\phi_k(a_{0}\dots a_{k}x_e) }\right)\times \left(\sum_{n=0}^{\infty}\sum_{\substack{(a_0,\dots,a_n)\in \mathcal{W}^{n+1}\\a_0=e,a_n=b\\
\forall i\geq 1:a_i\not\in A  }}\lambda^{-n}e^{\phi_n(a_{0}\dots a_{n}x_b) }\right)\\
 =&\sum_{e\in A} G(a,e|\lambda)L^A(e,b|\lambda).
\end{align*}
Moreover, all estimates are uniform in $a$ and $b$.
\item
Let $A \subseteq S$ be an arbitrary set. Then, 
\begin{align*}
& \sum_{e\in A}G(a,e|\lambda)L^A(e,b|\lambda)\\
=& \sum_{e\in A} \left(\sum_{k=0}^\infty \sum_{\substack{(a_0, \dots, a_k)\in \mathcal{W}^{k+1}\\ a_0=a, a_k = e}}\lambda^{-k}e^{\phi_k(a_0 \dots a_kx_{e})} \right)\times  \left(\sum_{n=0}^\infty \sum_{\substack{(a_0, \dots, a_n)\in \mathcal{W}^{n+1}\\ a_0=e, a_n = b\\  \;\forall i\geq 1:a_i\not\in A}}\lambda^{-n}e^{\phi_n(a_0 \dots a_nx_{b})} \right)\nonumber\\
\asymp& \sum_{e\in A}\sum_{n=0}^{\infty}\sum_{k=0}^n\sum_{\substack{(a_0, \dots, a_n)\in \mathcal{W}^{n+1}\\ a_0=a, a_k = e, a_{n}=b\\ \;\forall i> k:a_i\not\in A}}\lambda^{-n}e^{\phi_n(a_0 \dots a_nx_{b})}\nonumber \\
=& \sum_{n=0}^{\infty}\sum_{\substack{(a_0, \dots, a_n)\in \mathcal{W}^{n+1}\\ a_0=a, a_{n}=b\\ \;\exists i : a_i\in A}}\lambda^{-n}e^{\phi_n(a_0 \dots a_nx_{b})}.\nonumber
\end{align*}
Similarly, 
\begin{align}
\label{eq:FA_G}
& \sum_{e\in A}F^A(a,e|\lambda)G(e,b|\lambda)\\
=& \sum_{e\in A} \left(\sum_{k=0}^\infty \sum_{\substack{(a_0, \dots, a_k)\in \mathcal{W}^{k+1}\\ a_0=a, a_k = e\\ \forall i<k:a_i\not\in A}}\lambda^{-k}e^{\phi_k(a_0 \dots a_kx_{e})} \right)\times  \left(\sum_{n=0}^\infty \sum_{\substack{(a_0, \dots, a_n)\in \mathcal{W}^{n+1}\\ a_0=e, a_n = b\\ }}\lambda^{-n}e^{\phi_n(a_0 \dots a_nx_{b})} \right)\nonumber\\
\asymp&\sum_{e\in A}\sum_{n=0}^{\infty}\sum_{k=0}^n\sum_{\substack{(a_0, \dots, a_n)\in \mathcal{W}^{n+1}\\ a_0=a, a_k = e, a_{n}=b\\ \;\forall i< k:a_i\not\in A}}\lambda^{-n}e^{\phi_n(a_0 \dots a_nx_{b})}\nonumber \\
=&\sum_{n=0}^{\infty}\sum_{\substack{(a_0, \dots, a_n)\in \mathcal{W}^{n+1}\\ a_0=a, a_{n}=b\\ \;\exists i : a_i\in A}}\lambda^{-n}e^{\phi_n(a_0 \dots a_nx_{b})}\nonumber.
\end{align}
Therefore, 
$$ \sum_{e\in A}G(a,e|\lambda)L^A(e,b|\lambda) \asymp \sum_{e\in A}F^A(a,e|\lambda)G(e,b|\lambda).$$
Moreover, all estimates are uniform in $a,b$ and $A$.
\item
Let $A \subseteq S$ be an arbitrary set. Then, 
$$G(a,b|\lambda) \geq \sum_{n=0}^{\infty}\sum_{\substack{(a_0,\dots,a_n)\in \mathcal{W}^{n+1}\\
a_0=a, a_n=b\\\exists i: a_i \in A }}\lambda^{-n}e^{\phi_n(a_{0}\dots a_{n}x_b) }. $$
and, by Eq. (\ref{eq:FA_G}),
\begin{align*}
G(a,b|\lambda) 
\gg \sum_{e\in A}F^A(a,e|\lambda)G(e,b|\lambda).
\end{align*}
Moreover, all estimates are uniform in $a,b$ and $A$.
\item
We study the Green's function as a linear operator on  $C_c(X)$. 
Let $\T(a,b):
C_c(X)\rightarrow C_c(X)$, 
$$\T(a,b)f := 1_{[b]}L_\phi(1_{[a]}f) $$ 
and let $\T = \bigl(\T(a,b)\bigr)_{a,b\in S}$ be a $S\times S$ (infinite) matrix of operators. Let $\I$\ be the identity matrix, namely for every $a,b\in S$ and $f\in C_c(X)$, 
$$\I(a,b)f=\begin{cases}f & a=b \\
0 & a\neq b \\
\end{cases}.$$
Given two $S\times S$ matrices of operators ${\bf A}$ and ${\bf B}$, we define their $\star$-product by 
$${\bf (A} \star{\bf B)}(a,b) := \sum_{c\in S}{\bf B}(c,b){\bf A}(a,c). $$
Easy to verify that this product is associative. Then, for all $n>0$ and $f\in C_c(X)$, 
\begin{align*}
\T^{(n)} (a,b)f =& \underbrace{(\T \star \cdots \star \T)}_{n \text{ times}}(a,b)\\
=&\sum_{\substack{(a_0,\dots, a_n)\in \mathcal{W}^{n+1}\\ a_0=a, a_{n}=b }}\bigl(\T(a_{n-1}, a_n)\cdots \T(a_0, a_1)\bigr)f \\
=&1_{[b]}L_\phi^{n}(1_{[a]}f). 
\end{align*}
We let $\T^{(0) } = \I$.   
Let
$$\G_\lambda := \sum_{n=0}^\infty \lambda^{-n}\T^{(n)}. $$
Notice that 
$G(a,b|\lambda) =     (\G_{\lambda}(a,b)1)(bx_b)$.  For every $\lambda > \rho(\phi)$, 
\begin{align*}
\G_t \star \bigl(\I-\lambda^{-1}\T\bigr)=&\left( \sum_{n=0}^\infty  \lambda^{-n}\T^{(n)}\right)\star \bigl( \I - \lambda^{-1} \T \bigr)\\
 =& \sum_{n=0}^\infty \lambda^{-n}\T^{(n)} - \sum_{n=0}^\infty \lambda^{-n-1} \T^{(n+1)}\nonumber = \I.
\end{align*}
Similarly, $ \bigl(\I-\lambda^{-1}\T\bigr)\star \G\lambda = \I$. 
Therefore, for all $\lambda_1,\lambda_2 > \rho(\phi)$,
\begin{align}
\label{eq:G_t_almost}
&(\I - \lambda_1^{-1}\T)\left( \frac{\G_{\lambda_1}}{\lambda_1} - \frac{\G_{\lambda_2}}{\lambda_2} \right) (\I - \lambda_2^{-1} \T)\\
=&\lambda_1^{-1}\I-\lambda_1^{-1}\lambda_2^{-1} \T - \lambda_2^{-1} \I+ \lambda_{1}^{-1}\lambda_{2}^{-1}\T \nonumber\\
=&\left(\frac{1}{\lambda_1} - \frac{1}{\lambda_2}\right) \I.\nonumber
\end{align}
 We apply   $\G_{\lambda_1}$ on the left and  $\G_{\lambda_2}$ on the right to Eq. (\ref{eq:G_t_almost}) to obtain that
\begin{equation*}
  \frac{\G_{\lambda_1}}{\lambda_1} - \frac{\G_{\lambda_2}}{\lambda_2} = \left(\frac{1}{\lambda_1} - \frac{1}{\lambda_2}\right)(\G_{t_1}\star \G_{t_2} )
\end{equation*}
namely for every $a,b\in S$,
\begin{equation}
\label{eq:Gt_equation}
  \frac{\G_{\lambda_1}(a,b)}{\lambda_1} - \frac{\G_{\lambda_2}(a,b)}{\lambda_2} = \left(\frac{1}{\lambda_1} - \frac{1}{\lambda_2}\right)\sum_{c\in S}\G_{\lambda_2}(c,b) \G_{\lambda_1}(a,c).
\end{equation}
Notice that
$$\bigl(\T^{(m)}(c,b)\T^{(n)}(a,c)\bigr)1(bx_b)=L_\phi^{m}\bigl(1_{[c]}L_\phi^{n}(1_{[a]})\bigr)(bx_{b})\asymp  L_\phi^{m}(1_{[c]})(bx_b)L_\phi^{n}(1_{[a]})(cx_c).$$ 
Then, the proposition follows by  Eq. (\ref{eq:Gt_equation}) evaluated on the function $f \equiv 1$ and at the point $bx_b$. 
\item
Since $b_i\neq a_j$ for  $i=1,2$ and $1\leq j\leq N$,
\begin{align*}
&G(1_{[a_1,\dots, a_N]}, b_i|\lambda) =\sum_{n=N}^\infty\sum_{\substack{(c_{0},\dots,c_n)\in \mathcal{W}^{n+1}\\ c_n=b_{i}\\\forall k<N:c_k = a_{k+1}}}\lambda^{-n}e^{\phi_n(c_{0},\dots,c_nx_{b_i})}\\
\asymp& \sum_{n=N}^\infty\sum_{\substack{(c_{0},\dots,c_n)\in \mathcal{W}^{n+1}\\ c_n=b_{i}\\\forall k<N:c_k = a_{k+1}}}\lambda^{-(N-1)}e^{\phi_{N-1}(c_{0},\dots,c_Nx_{c_{N}})}\lambda^{-n-N+1}e^{\phi_{n-N+1}(c_{N},\dots,c_nx_{b_i})}\\
=&\ \lambda^{-(N-1)}e^{\phi_{N-1}(a_1\dots a_N x_{a_N})}\sum_{n=0}^\infty\sum_{\substack{(c_{0},\dots,c_n)\in \mathcal{W}^{n+1}\\ c_{0 }= a_N, c_n=b_{i}\\}}\lambda^{-n}e^{\phi_{n}(c_{0},\dots,c_nx_{b_i})}\\
=& \ \lambda^{-(N-1)}e^{\phi_{N-1}(a_1\dots a_N x_{a_N})} G(a_N, b_i|\lambda).
\end{align*}
Hence, 
$$K(1_{[a_1,\dots, a_N]}, b_i|\lambda)\asymp \lambda^{-N-1}e^{\phi_{N-1}(a_1\dots a_N x_{a_N})}K(a_N, b_i|\lambda) $$
and by Proposition \ref{prop:green_main_prop}.\ref{prop:green_main_prop_FG},
\begin{align*}
\frac{K(1_{[a_1,\dots, a_N]}, b_1|\lambda)}{K(1_{[a_1,\dots, a_N]}, b_2|\lambda)}\asymp  \frac{K(a_N, b_1|\lambda)}{K(a_N, b_2|\lambda)}\asymp \frac{F(a_N, b_1|t)F(o, b_2|\lambda)}{F(o, b_1|t)F(a_{N}, b_2|\lambda)}.
\end{align*}
Moreover, all estimates are uniform in $b_1, b_2$ and in $[a_1,\dots, a_N]$.
\end{enumerate} 
\end{proof}
\begin{prop}
\label{prop:green_main_prop_2}
Assume that $(X,T)$ is locally compact and transitive, that $\phi$ has summable variations, that  $\phi$ is uniformly irreducible w.r.t. a connected, undirected and locally finite graph $(S,E)$ and  that $P_G(\phi)<0$. Then,
\begin{enumerate}[label={(\arabic*)}]
\item 
\label{prop:green_main_prop_L}
For every $a,b\in S$ and every $\lambda\in \bigl(\rho(\phi), 1\bigr)$, 
$$L^A(a,b|\lambda) \geq \lambda^{-d_{E}(a,b)}L^A(a,b).$$
\item
\label{prop:green_main_prop_harnack}
(Harnack's inequality) There exists $ C' >1$ s.t. 
for every $\lambda \in \bigl(\rho(\phi),1\bigr]$,  for every $h\in \bigl\{G(\cdot, c|\lambda),F(\cdot,c|\lambda),G(c, \cdot|\lambda),F(c, \cdot|\lambda)\bigr\}_{c\in S}$ and  every $a,b\in S$,
$$h(a) \leq ( C')^{  d_{E}(a,b)}h(b). $$
\end{enumerate}
\end{prop}
\begin{proof}
\begin{enumerate}[label={(\arabic*)}]
\item
Let $\lambda\in \bigl( \rho(\phi),1\bigr) $ and let $k$ the minimal number s.t. $L_\phi^k(1_{[a]})(bx_b)>0$. Since  $\phi$ is uniformly-irreducible w.r.t. the set of edges $E$, $k\geq d_E(a,b)$ and 
\begin{align*}
L^A(a,b|\lambda)=&\sum_{n=0}^\infty\sum_{\substack{(a_0, \dots, a_n)\in \mathcal{W}^{n+1}\\ a_0=a, a_n = b\\a_0 \in A, \;\forall i\geq 1:a_i\not\in A }}\lambda^{-n}e^{\phi_n(a_0 \dots a_nx_{b})}\\
=& \sum_{n=d_{E}(a,b)}^\infty\sum_{\substack{(a_0, \dots, a_n)\in \mathcal{W}^{n+1}\\ a_0=a, a_n = b\\a_0 \in A, \;\forall i\geq 1,a_i\not\in A }}\lambda^{-n}e^{\phi_n(a_0 \dots a_nx_{b})}\\
\geq&\lambda^{-d_E(a,b)} \sum_{n=d_{E}(a,b)}^\infty\sum_{\substack{(a_0, \dots, a_n)\in \mathcal{W}^{n+1}\\ a_0=a, a_n = b\\a_0 \in A, \;\forall i\geq 1,a_i\not\in A }}e^{\phi_n(a_0 \dots a_nx_{b})}\\
=&\lambda^{-d_E(a,b)} L^{A}(a,b).
\end{align*}
\item
Let $K>0$ s.t. for every $(a,b)\in E$ there exists $k\leq K$ with$$L_\phi^k(1_{[a]})(bx_b)>0 $$
and let  
$$\epsilon = \exp\bigl(-K\min_{x\in X}|\phi(x)|\bigr). $$
Let $N=d_E(a,b)+1$ and let $a_1,\dots, a_N$ be a shortest path  in $E$ from $a_1=b$ to $a_N=a$. Let $k_1,\dots, k_{N-1}\geq 1  $ with $k_i\leq K$ and $L_\phi^{k_i}(1_{[a_i]})(a_{i+1} x_{a_{i+1}})>0$. Observe that for every $i$, 
$$ L_\phi^{k_i}(1_{[a_i]})(a_{i+1} x_{a_{i+1}}) \geq \epsilon.   $$
 Let $k = \sum_{i=1}^{N-1}k_i $. Notice that $d_E(a,b) \leq k\leq d_E(a,b)K$. Then, \begin{align*}G(b,c|\lambda) \geq& \sum_{n=k}^\infty \lambda^{-n}L_\phi^{n}(1_{[b]})(c x_c) \\
\geq &\ \lambda^{-k}\sum_{n=0}^\infty \lambda^{-n}L_\phi^{n+k}(1_{[b]}\cdot 1_{T^{-d_E(a,b)}[a]})(c x_c)\\
\geq & \lambda^{-k}C_\phi^{-k}\sum_{n=0}^\infty L_\phi^{k_{1}}(1_{[b]})(a_{2} x_{a_2})\cdots L_\phi^{k_{N-1}}(1_{[a_{N-1}]})(a x_{a})\lambda^{-n} L_\phi^{n}(1_{[a]})(c x_c)\\
\geq &\ \lambda^{-k}C_\phi^{-k}\epsilon_0^{d_E(a,b)}G(a,c|\lambda)\\
\geq &\ \bigl( \max\{C_\phi, C_\phi^{K}\} \epsilon^{-1} \bigr)^{-d_E(a,b)}.  
\end{align*}
So, with $   B= \max\{C_\phi, C_\phi^{K}\} \epsilon_0^{-1}$,
$$G(a,c|\lambda) \leq B ^{d_E(a,b)}G(b,c|\lambda). $$
Similar arguments lead to the following inequality$$G(c,b|\lambda) \leq B^{d_E(a,b)}G(c,a|\lambda) . $$ 
 Let $C>1$ be the constant from Proposition \ref{prop:green_main_prop}.\ref{prop:green_main_prop_FG}. Then, 
\begin{align*}
F(a,c|\lambda) \geq& C^{-1} \frac{G(a,c|\lambda)}{G(c,c|\lambda)} \\
\geq & C^{-1}B^{-d_E(a,b)}\frac{G(b,c|\lambda)}{G(c,c|\lambda)}\\
\geq & C^{-2}B^{-d_E(a,b)}F(b,c|\lambda).
\end{align*}
Moreover, since
$$G(a,a|\lambda) \geq B^{-d_E(a,b)}G(a,b|\lambda)  $$
and
$$G(b,b|\lambda) \leq B^{d_E(a,b)}G(a,b|\lambda) $$
we have that
$$\frac{G(a,a|\lambda)}{G(b,b|\lambda)}   \geq  B^{-2d_E(a,b)}.$$
We conclude, 
\begin{align*}
F(c,b|\lambda) \geq& C^{-1} \frac{G(c,b|\lambda)}{G(b,b|\lambda)} \\
\geq & C^{-1}B^{-d_E(a,b)}\frac{G(c,a|\lambda)}{G(b,b|\lambda)}\\
\geq & C^{-2}B^{-d_E(a,b)}F(c,a|\lambda)\frac{G(a,a|\lambda)}{G(b,b|\lambda)}\\
\geq &\ C^{-2}B^{-3d_E(a,b)}F(c,a|\lambda). 
\end{align*}
\end{enumerate}
\end{proof}
\subsection{Proof of Theorem \ref{thm:hyperbolic_main_thm}}
 We follow here the arguments of the proof of the original theorem as presented in \cite{woess_2000}.
  
Recall that if $\lambda > \rho(\phi)$ then $P_G(\phi - \log \lambda) <\ 0$ and that for
 all $f\in C_c(X)$ and $x\in X$,
$$\lambda^{-n} (L_\phi^n f)(x) = L_{\phi - \log \lambda}^n (f) (x). $$
Thus we
can assume w.l.o.g. that $P_G(\phi) < 0$ and prove the theorem for $\lambda=1$. 

In what follows, assume that $(S,E)$ is a  $\delta$-hyperbolic graph and that $P_G(\phi)<0$.  For $a,b\in S$, let
$$U_{a,b} = \bigl\{c\in S : |b\wedge c|_{a} \geq d_{E}(a,b)-7\delta\bigr\} $$
and let $V_{b,a} = S \setminus U_{a,b}$. For $a\in S$ and $r\geq 0$, we denote by 
$ B(a,r) = \{b\in S: d_{E}(a,b)\leq r\}$ the closed ball of radius $r$ around $a$. Let $C,C'>1$ be the constants from Propositions \ref{prop:green_main_prop} and \ref{prop:green_main_prop_2} respectively and let $C_0 =\max \{C, C'\}$.
\begin{prop}
\label{prop:UV_ineq}
Under the assumptions of Theorem \ref{thm:hyperbolic_main_thm}, for every   $\lambda\in \bigl(\rho(\phi), 1\bigr)$  there exists a constant $C_{1}(\lambda)>1$ s.t. for every $a,b\in S$ and for every $v$ on some geodesic segment from $a$ to $b$,
$$ G(a,w)\leq C_{1}(t)F(a,v)G(v,w|\lambda) , \quad \forall w\in U_{a,v}\cup V_{v,b}.$$
\end{prop}
\begin{proof}
We use only properties that do not depend on the base point and so we can assume w.l.o.g. that $a=o$. Let $\ell=21\delta$, let $m$ be the integer part of $d(o,v)/\ell$ ($m$ may be zero) and consider the points $v_0,\dots, v_m$ which lie on a geodesic segment between $o$ and $v$ with $d_{E}(v_k,v) = (m-k)\ell$. Let $W_k = U_{o,v_k}\cup V_{v_k, b}$ and let $dW_k = \{w\in W_k:d_{E}(w, S\setminus W_k) = 1\}$. 
\begin{lemma}
\label{lemma:woess_Wk}
Assume that $(S, E)$ is a    $\delta$-hyperbolic graph. Then, for all $k\geq 1$, 
\begin{enumerate}
\item 
$v_k\in  W_k \subseteq W_{k-1}$. 
\item
If $w\in W_k$ with $d_{E}(w, v_k)\geq 2r+\ell+1$ then $B(w,r) \subseteq W_{k-1}$.
\end{enumerate}
\end{lemma}
 \begin{proof}
 See Lemma 27.7 in \cite{woess_2000}.
 \end{proof}

Choose an integer $r\geq \ell$ with $\lambda^{r}C_0^{2l+4} \leq 1$. We show by induction on $k$ that, with $C_1 =C_0^{4r+2\ell+1}  $, 
\begin{equation}
\label{eq:Wk_ineq}
G(o,w) \leq C_{1}F(o,v_k) G(v_k,w|\lambda), \quad \forall w\in W_k. \end{equation}
The proposition follows with $k=m$.

 Let $k=0$. Then, $d_{E}(o, v_0)\leq \ell$ and by Proposition \ref{prop:green_main_prop_2}.\ref{prop:green_main_prop_harnack}, for every $w\in S$, 
$$G(o,w)\leq C_0^\ell G(v_0, w)\leq C_0^\ell G(v_0, w|\lambda). $$
Similarly, $G(v_0, v_0)\leq C_0^\ell G(o,v_0)$. By Proposition \ref{prop:green_main_prop}.\ref{prop:green_main_prop_FG}, \begin{equation}
\label{eq:F_geq_1}
 F(o, v_0) \geq C_0^{-1} \frac{G(o,v_0)}{G(v_0,v_0)}\geq C_0^{-\ell-1}
\end{equation}
and thus, for every $w\in S$,
$$ G(o,w)\leq C_0^{2\ell+1} F(o, v_0)G(v_0, w|\lambda).$$

Next, suppose by induction that Eq. (\ref{eq:Wk_ineq}) holds for $k-1$. Since $d_{E}(v_{k-1}, v_k)\leq \ell$, by Proposition \ref{prop:green_main_prop_2}.\ref{prop:green_main_prop_harnack}, for all $w\in S$,
\begin{equation}
\label{eq:G_ineq_reduction}
G(v_{k-1}, w|\lambda)\leq C_0^{\ell}G(v_k,w|\lambda).
\end{equation}
Similarly to Eq. (\ref{eq:F_geq_1}),
\begin{equation}
\label{eq:F_geq_1_2}
C_0^{\ell+1}F(v_{k-1}, v_k)\geq1. 
\end{equation}
Hence,  for all $w\in W_{k-1}$,
\begin{align}
\label{eq:W_(k-1)induction_ineq}
\\
G(o,w)\leq
& C_{1}F(o, v_{k-1})G(v_{k-1}, w|\lambda) & (\because \; \text{induction hypothesis)}\nonumber\\
\leq  & C_{1}C_0^{2\ell+1} F(o, v_{k-1})F(v_{k-1}, v_k)G(v_{k}, w|\lambda) & (\because \; \text{Eq. }(\ref{eq:G_ineq_reduction}, \ref{eq:F_geq_1_2}))
 \nonumber\\
\leq  & C_{1}C_0^{2\ell+2} F(o,  v_k)G(v_{k}, w|\lambda). & (\because \; \text{Proposition }\ref{prop:green_main_prop}.\ref{prop:green_main_prop_F}) \nonumber
\end{align}
Now, let $w\in W_k$ and assume first that $d_{E}(w,v_k)\geq 2r+\ell+1$. Set $A=\{e\in S: d_E(e,w)=r\}$.  By Lemma \ref{lemma:woess_Wk}, $A \subseteq W_{k-1}$ and hence  Eq. (\ref{eq:W_(k-1)induction_ineq}) holds for all $e\in A$. 
We claim that any path from $o$ to $w$ must pass through $A$. If $k=1$ then, by construction, 
$$d_E(o,v_1)\geq \ell $$
and, by assumption,  
$$d_E(w, v_1)\geq 2r+\ell+1 $$
so $d_E(o,w) >r$. Thus any path from $o$ to $w$ must enter $A$. Observe that $o\not\in W_1$\ and by Lemma  \ref{lemma:woess_Wk} $o\not\in W_k$ for all $k$. Hence, if $k\geq 2$ then $o \not \in B(w,r)$ and again any path from $o$ from $w$ must enter $A$.
  
We deduce, 
\begin{align*}G(o,w) \leq& C_0\sum_{e\in A}G(o,e) L^A(e,w)& (\because \; \text{Proposition }\ref{prop:green_main_prop}.\ref{prop:green_main_prop_GL})
\\
\leq & C_{1}C_0^{2\ell+3} \sum_{e\in A} F(o,  v_k)G(v_{k}, e|\lambda) L^A(e,w) &  (\because \;\text{Eq. }(\ref{eq:W_(k-1)induction_ineq})) 
\\
\leq & C_{1}C_0^{2\ell+3}\lambda^{r} \sum_{e\in A} F(o,  v_k)G(v_{k}, e|\lambda) L^A(e,w|\lambda)
.& (\because \; \text{Proposition }\ref{prop:green_main_prop_2}.\ref{prop:green_main_prop_L})
 \end{align*}
 Since $d(w,v_k)>r$, any path from $w$ to $v_k$ must pass through $A$. Therefore,
by Proposition \ref{prop:green_main_prop}.\ref{prop:green_main_prop_GL}, $$C_{1}C_0^{2\ell+3}\lambda^{r} \sum_{e\in A} F(o,  v_k)G(v_{k}, e|\lambda) L^A(e,w|\lambda)\leq C_{1}C_0^{2\ell+4}t^{r}  F(o,  v_k)G(v_{k}, w|\lambda). $$
So, for all $w\in W_{k}$ with $d_{E}(w,v_k)\geq 2r+\ell+1$,
$$ G(o,w)\leq C_{1}C_0^{2\ell+4}\lambda^{r}  F(o,  v_k)G(v_{k}, w|\lambda).  $$
By the choice of $r$, Eq. (\ref{eq:Wk_ineq}) follows.

Lastly, if $w\in W_k$ with $d_{E}(w,v_k)\leq 2r+\ell$, then\begin{align*}
G(o,w) \leq& C_0^{2r+\ell}G(o,v_k)&
 (\because \; \text{Proposition }\ref{prop:green_main_prop_2}.\ref{prop:green_main_prop_harnack}) \\
\leq &  C_0^{2r+\ell+1} F(o,v_k)G(v_{k}, v_k)
&  (\because \; \text{Proposition }\ref{prop:green_main_prop}.\ref{prop:green_main_prop_FG})\\
\leq &  C_0^{2r+\ell+1} F(o,v_k)G(v_{k}, v_k|\lambda)
&  (\because \; \lambda<1)\\
\leq & C_0^{4r+2\ell+1}F(o,v_k)G(v_{k}, w|\lambda). 
&  (\because \; \text{Proposition }\ref{prop:green_main_prop_2}.\ref{prop:green_main_prop_harnack})
\end{align*} 
\end{proof}

\begin{cor}
\label{cor:VU_ineq}
Under the assumptions of Theorem \ref{thm:hyperbolic_main_thm}, for every   $\lambda\in \bigl(\rho(\phi), 1\bigr)$  there is a constant $C_{2}(\lambda)>1$ s.t. for every $a,b\in S$ and for every $v$ on some geodesic segment from $a$ to $b$,
$$ G(w,b)\leq C_{2}(\lambda)G(w,v|\lambda)L(v,b) , \quad \forall w\in V_{v,a}\cup U_{b,v}.  $$
\end{cor}
\begin{proof}
Denote by  $X^{\pm}$ is the \textit{two-sided shift}, by $X^-$  the \textit{negative one-sided shift} and by $X^+=X$  the \textit{positive one-sided shift}. It is known that there exists a potential function $\phi^-:X^{-} \rightarrow\mathbb{R}$  with summable variations and a bounded uniformly continuous function $\psi:X^{\pm}\rightarrow\mathbb{R}$ s.t.
$$\phi^+ - \phi^- = \psi - \psi\circ T. $$   
 See Section 5 in \cite{shwartz_2019}. We add the $+$ or $-$ notation over the Green's functions to distinct between the two spaces, e.g. $G^+$ or $G^-$.
\begin{lemma}
\label{lemma:G_plus_minus}
Assume that $(X,T)$ is locally compact and transitive, that $\phi$ has summable variations and that $P_G(\phi)<\infty$. Then, 
 there exists $C_{}''>1$ s.t. for every $\lambda>\rho(\phi)$ and every $a,b,c\in S$,
 $$G^-(a,b|\lambda) =(C'')^{\pm 1} G^+(b,a|\lambda)$$
and 
$$F^{-}(a,b) G^-(b,c|\lambda) =(C'')^{\pm 1} G^+(c,b|\lambda) L^+(b,a|\lambda). $$
\end{lemma}
\begin{proof}
Let $a,b\in S$, let $a_1,\dots, a_{n-1}$ be an admissible path from $a$ to $b$, let $x_a^-\in [a] \subseteq X^-$ and let $x_b^+ \in [b] \subseteq X^+$. Since $\phi^-, \phi^+$ and $\psi$ are all bounded, \begin{align*}
&|\phi_n^+(a, a_1, \dots, a_{n-1}x_b^+) - \phi_n^-(x_a^- a_1, \dots, a_{n-1},b)|  \\
&\leq |\phi^+(a ,a_1, \dots, a_{n-1}x_b^+)  |+ |\phi^-(x_a^- a_1, \dots, a_{n-1},b)|\\
&+\left| \sum_{i=1}^{n-1}\left( \phi^+ (a_i, \dots, a_{n-1}x_b^+ )-\phi^- (x_a^-a_1,\dots, a_i)\right) \right|\\
&\leq  \sup|\phi^+| + \sup|\phi^-| +2 \sup|\psi|. 
\end{align*}
Now, any path from $a$ to $b$ in $X^+$ is a path from $b$ to $a$ in $X^-$ and there thus is a natural matching of the terms in the sums $G^+, G^-$ and $F^+, F^-$ with the property that matching terms are within multiplicative error $e^{\pm (\sup|\phi^+| + \sup|\phi^-| +2 \sup|\psi|)}$ from each other.   
\end{proof}
The corollary follows from Proposition \ref{prop:UV_ineq} and Lemma \ref{lemma:G_plus_minus}.
 \end{proof}

\begin{thm}
\label{thm:main_F_ineq}(Ancona's inequality)
Under the assumptions of Theorem \ref{thm:hyperbolic_main_thm}, for every $r\geq 0$ there exists $C_{3}(r)\geq 1$ s.t. 
$$\bigl(C_{3}(r)\bigr)^{-1}F(a,v)F(v,b) \leq F(a,b) \leq C_{3} (r)F(a,v)F(v,b) $$
whenever $a,b\in S$ and $v$ is at distance at most $r$ from some geodesic segment from $a$ to $b$. 
\end{thm}
\begin{proof}
The lower bound follows from Proposition \ref{prop:green_main_prop}.\ref{prop:green_main_prop_F}, so we focus on the upper bound. We first consider the case $r=0$, when $v$ lies on a geodesic segment from $a$ to $b$. 

Fix $\lambda\in \bigl(\frac{\rho(\phi)+1}{2}, 1\bigr)$.
 If $d_{E}(a,v)\leq 7\delta$ then by Proposition \ref{prop:green_main_prop_2}.\ref{prop:green_main_prop_harnack}, $$F(a,b)\leq  C_0^{7\delta}F(v,b)$$ 
and
$$1\leq F(v,v)\leq C_0^{7\delta} F(a,v).$$
In particular,  
$$F(a,b)\leq C_0^{14\delta}F(a,v)F(v,b). $$

Suppose that $d_{E}(a,v)> 7\delta$. Since $|a\wedge v|_a=0$, $a\not\in U_{a,v}$. Moreover, since $v$ lies on a geodesic segment from $a$ to $b$, $|v\wedge b|_a = d_E(a,v)$ and thus $b\in U_{a,v}$. In particular,  any path from $a$ to $b$ must pass through 
$$A:= \bigl\{c\in U_{a,v} :\exists w \in S \setminus U_{a,v}, \; d_E(c,w) = 1\bigr\}. $$
 By Propositions \ref{prop:green_main_prop}.\ref{prop:green_main_prop_GL} and \ref{prop:UV_ineq},
\begin{equation}
\label{eq:ancona_eq_1}
G(a,b) \leq C_0\sum_{w\in A}G(a,w)L^A(w,b) \leq C_0C_{1}F(a,v)\sum_{w\in A}G(v,w|\lambda)L^A (w,b) 
\end{equation}  
where $C_1 = C_1(1)$ is the constant from Proposition \ref{prop:UV_ineq}.
Every point $w\in A$ is at distance $1$ from some point $w'$ in $V_{v,a}=S \setminus U_{a,v}$. By Corollary \ref{cor:VU_ineq}, $$ G(w',b)\leq C_{2}G(w',v|\lambda) L(v,b) $$
where $C_2 = C_2(1)$ is the constant from Corollary 
\ref{cor:VU_ineq}. We apply Proposition \ref{prop:green_main_prop_2}.\ref{prop:green_main_prop_harnack} to $G(w', b)$ and $G(w', v|\lambda)$  to obtain that
\begin{equation}
\label{eq:harnack_GGL_ineq}
G(w,b)\leq C_{2}C_0^2 G(w,v|\lambda)L(v,b). 
\end{equation}
Then, by Proposition \ref{prop:green_main_prop}, for every $e\in S$, 
\begin{align}
\label{eq:ancona_eq_2}
\\
\sum_{w\in A}G(e,w)L^A(w,b)\leq & C_0 \sum_{w\in A}F^A(e,w)G(w,b)&
(\because \; \text{Proposition }\ref{prop:green_main_prop}.\ref{prop:green_main_prop_GLFG})
\nonumber\\
\leq  & C_0 \sum_{w\in A}F^A(e,w|\lambda)G(w,b)&
(\because \; \lambda<1)\nonumber\\
\leq  & C_{2} C_0^3  \sum_{w\in A}F^A(e,w|\lambda)G(w,v|\lambda)L(v,b)&
(\because \; \text{Eq. }(\ref{eq:harnack_GGL_ineq}))\nonumber\\
\leq & C_{2} C_0^4 G(e,v|\lambda)L(v,b).&
(\because \; \text{Proposition }\ref{prop:green_main_prop}.\ref{prop:green_main_prop_FEG})
\nonumber
\end{align}
Let $\nu(e) = \lambda \delta_v(e) + (1-\lambda) G(v, e|\lambda)$. By Proposition \ref{prop:green_main_prop}.\ref{prop:green_main_prop_Gt}, with $\lambda_1=1$ and $\lambda_2=\lambda$
\begin{align}
\label{eq:ancona_eq_3}
G(v,w|\lambda) \leq &\lambda G(v,w) + C_0(1-\lambda)\sum_{e\in S} G(v,e|\lambda)G(e,w)\nonumber\\
\leq &C_0 \left ( \lambda G(v,w) + (1-\lambda)\sum_{e\in S} G(v,e|\lambda)G(e,w) \right)\\
=& C_{ 0 }\sum_{e\in S}\nu(e) G(e,w).\nonumber
\end{align}
In summary, 
\begin{align*}
G(a,b) \leq &C_0C_{1}F(a,v)\sum_{w\in A}G(v,w|\lambda)L^A (w,b) &(\because \; \text{Eq. }(\ref{eq:ancona_eq_1}))\\
\leq & C_{0}^{2}C_{1}F(a,v)\sum_{w\in A}\sum_{e\in S}\nu(e) G(e,w)L^A (w,b)
&(\because \; \text{Eq. }(\ref{eq:ancona_eq_3}))\\
= & C_{0}^{2}C_{1}F(a,v)\sum_{e\in S}\nu(e)\left(\sum_{w\in A} G(e,w)L^A (w,b)\right)\\
\leq & C_{0}^{6}C_{1}C_{2}F(a,v)\left(\sum_{e\in S}\nu(e)G(e,v|\lambda)\right)L(v,b).
& (\because \; \text{Eq. }(\ref{eq:ancona_eq_2}))
\end{align*}
Choose $\lambda_2(\lambda) \in\bigl( \rho(\phi), 2\lambda-1\bigr)$ with $\lim_{\lambda\rightarrow 1^-}\lambda_2(\lambda)=1$.   Then, 
$$1-\lambda \leq \lambda - \lambda_2 $$
and
\begin{align*}
\sum_{e\in S}\nu(e) G(e, v|\lambda) = & \lambda G(v,v|\lambda) + \left(1-\lambda\right)\sum_{e\in S}G(v,e|\lambda)G(e,v|\lambda)&\\
\leq&\lambda G(v,v|\lambda) +(\lambda-\lambda_{2})\sum_{e\in S}G(v,e|\lambda_{2})G(e,v|\lambda) &\\
\leq & \lambda G(v,v|\lambda) + C_0 \frac{\lambda-\lambda_2}{\frac{1}{\lambda_2}-\frac{1}{\lambda}}\left(\frac{G(v,v|\lambda_2)}{\lambda_2}-\frac{G(v,v|\lambda)}{\lambda}
\right)& (\because \; \text{Proposition }\ref{prop:green_main_prop}.\ref{prop:green_main_prop_Gt})\\
\leq & \lambda G(v,v|\lambda) + C_0t G(v,v|\lambda_{2})\\
\leq &\ C_0 t\bigl(G(v,v|\lambda) + G(v,v|\lambda_{2})\bigr).
\end{align*}
This leads to $$G(a,b) \leq  C_{0}^{7}C_{1}C_{2}tF(a,v)\bigl(G(v,v|\lambda) + G(v,v|\lambda_{2})\bigr)L(v,b). $$
Since $\rho(\phi)<1$,  $G(v,v|\lambda)$ is analytic as a function of $\lambda$ on a neighbourhood of $\lambda=1$ and $$\lim_{\lambda \nearrow 1} G(v,v|\lambda)=\lim_{\lambda \nearrow 1} G(v,v|\lambda_2(\lambda))=G(v,v).$$ Therefore$$G(a,b) \leq  2C_{0}^{7}C_{1}C_{2}F(a,v)G(v,v) L(v,b). $$
By Proposition \ref{prop:green_main_prop}.\ref{prop:green_main_prop_GL},
$$G(a,b)\leq 2C_0^7 C_1 C_2 F(a,v) G(v,b). $$
We divide both sides by $G(b,b)$ and apply Proposition \ref{prop:green_main_prop}.\ref{prop:green_main_prop_FG} to obtain that 
$$F(a,b) \leq 2C_0^9 C_1 C_2 F(a,v)F(v,b). $$
This proves Ancona's inequality in case $v$ lies on a geodesic segment from $a$ to $b$. 

Now, assume that $v$ is at distance $r\geq0$ from some geodesic segment from $a$ to $b$. Then, we can find $v'\in S$ on this geodesic segment from $a$ to $b$ with $d(v,v')=r$. By the first part of the proof,   
$$F(a,b)\leq 2C_0^9 C_1 C_2 F(a,v')F(v',b).  $$
Applying Proposition \ref{prop:green_main_prop_2}.\ref{prop:green_main_prop_harnack} twice leads to
$$F(a,b)\leq 2C_0^{9+2r}C_1 C_2 F(a,v)F(v,b).$$
\end{proof}

\paragraph{\bf{Proof of Theorem \ref{thm:hyperbolic_main_thm}.}}
Let $\xi \in \partial (S,E)$. 
 We first show that there exists $\epsilon_1\in (0,1)$ s.t. for every  $f\in C_c^+(X)$ and for every two sequences  $b_n, b_n'\in S$ which converge to  $\xi $,
 $$\liminf  _{n\rightarrow\infty}K(f, b_n)\geq \epsilon_1 \limsup_{n\rightarrow\infty} K(f, b'_n). $$
\begin{lemma} 
Let $(S, E)$ be a   $\delta$-hyperbolic graph. Then, for every $n$, there exists $v_{n}\in S$ which  is at distance at most $2\delta$ from some geodesic segments from $a$ to $b_n$, from $o$ to $b_n$, from $a$ to $b_n'$ and from $o$ to $b_n'$.
\end{lemma}
\begin{proof}
See  \cite{woess_2000}, proof of Theorem 27.1.
\end{proof}
For every $n>0$, let $v_n\in S$ as in the lemma and let $a_1,\dots, a_N\in S$ with $[a_1,\dots, a_N]\neq \varnothing$. By Proposition \ref{prop:green_main_prop}.\ref{prop:green_main_prop_K} and Theorem \ref{thm:main_F_ineq},
for all $n$ large enough,
\begin{align*}
\frac{K(1_{[a_1,\dots, a_N]}, b_n)}{ K(1_{[a_1,\dots, a_N]}, b_n')} \geq  &C_0^{-1} \frac{F(a_, b_n)F(o,b_n')}{F(a_, b_n')F(o, b_n)} \\\nonumber
\geq &C_0^{-1} \bigl(C_3(2\delta)\bigr)^{-4} \frac{F(a,v_{n})F(v_{n},b_n) F(o,v_{n})F(v_{n},b_n')}{F(a,v_{n})F(v_{n},b_n')F(o,v_{n}) F(v_{n}, b_n)} \\\nonumber
=& C_0^{-1} (C_3(2\delta))^{-4}.
\end{align*}    
Here $C_3(2\delta)$ is the constant in Ancona's inequality for $r=2\delta$. In particular,
\begin{equation}
\label{eq:hyperbolic_main_thm_ineq}
\liminf_{n\rightarrow\infty}  K(1_{[a_1,\dots, a_N]}, b_n) \geq  C_0^{-1} \bigl(C_3(2\delta)\bigr)^{-4}  \limsup_{n\rightarrow\infty} K(1_{[a_1,\dots, a_N]}, b'_n).
 \end{equation}
Since the collection cylinder sets linearly spans a dense subset of $C_c(X)$ w.r.t. the sup norm, Eq. (\ref{eq:hyperbolic_main_thm_ineq}) extends to all $f\in C_c^+(X)$. 

Let $\epsilon_1 =  C_0^{-1} \bigl(C_3(2\delta)\bigr)^{-4} \in (0,1)$ and let $$\mathcal{\mathcal{A}_\xi} = \bigl\{\omega  : \exists b_n\in S \text{ s.t. } b_n\xrightarrow[]{}\xi \text{}\text{ and } \lim_{n\rightarrow\infty}K(f,b_n) = K(f,\omega),\; \forall f\in C_c(X)\bigr\},$$ that is the set of all possible limit points in $\mathcal{M}$ of sequences  $b_n\in S$ with $b_n \rightarrow \xi$ in $(S,E)$. 
 We show that $\mathcal{A}_\xi$ consists of a single point alone. Notice that for all $\omega_1,\omega_2\in \mathcal{A}_\xi$, $$K(f,\omega_1) \geq \epsilon_1K(f, \omega_2), \quad \forall f\in C_c^+(X).$$ 
Therefore, it suffices to show that  $\mathcal{A}_\xi\cap \mathcal{M}_m\neq \varnothing$.

Let $\omega \in\mathcal{A}_\xi $, 
$$\mathcal{C} = \conf(1) =\left\{ \mu \text{ Radon}: \mu \geq 0 \text{ and } L_\phi^* \mu = \mu \right\} $$
and let 
 $$\mathcal{B}_\omega=\left\{\mu\in \mathcal{C} :   \sup_{ f\in C_c^+(X)}\frac{\mu(f)}{ \mu_{\omega}(f)} = 1\right\}.$$  Recall that $\mu_\omega(f) = K(f,\omega)$, $f\in C_c(X)$. If $\mu_{\omega} = \mu_1 + \mu_2$ with $\mu_i\in \mathcal{C}$ and the measures  $\mu_1,\mu_2$ are mutually singular and non-zero then  $\sup_{f\in C_c^+(X)}\bigl\{\frac{\mu_i}{\mu_{\omega}}\bigr\}= 1$. Thus, it suffices to show that $\mathcal{B}_\omega=\{\mu_\omega\}$.

Let $(b_0, b_1,\dots)$ be a geodesic sequence converging to $\xi$ in $(S,E)$ with $b_0=o$ and $\lim_{n\rightarrow\infty}K(f,b_n) = K(f,\omega)$ for every $f\in C_c(X)$. By Proposition \ref{prop:green_main_prop}.\ref{prop:green_main_prop_FG} and Theorem \ref{thm:main_F_ineq},
$$K(b_{k}, b_n)\geq C^{-2}_0 \frac{F(b_{k},b_n)}{F(o,b_n)}\geq C_0 ^{-2}C_3^{-1} \frac{1}{F(o,b_k)} .$$
where $C_3 = C_3(0)$.
Therefore, for every $a\in S$, 
$$K(a, b_k) \leq C^{2}_0 \frac{F(a, b_k)}{F(o, b_k)}\leq C_3^{}C_0^4F(a,b_k) K(b_{k}, b_{n})\xrightarrow[n\rightarrow\infty]{}C_3^{}C_0^4F(a,b_k) K(b_{k}, \omega). $$
 Let $a_1,\dots, a_N\in S$ admissible. By Proposition \ref{prop:green_main_prop}.\ref{prop:green_main_prop_K},  for all $k$ large enough, \begin{align}
\label{eq:omega_ext_eq}
 K(1_{[a_{1},\dots, a_N]}, b_k) \leq & C_0 e^{\phi_{N-1}(a_{1},\dots, a_{N}x_{a_N})}K(a_N, b_k)\\
 \leq &C^{5}_{0}C_3^{}e^{\phi_{N-1}(a_{1},\dots, a_{N}x_{a_N})}F(a_{N},b_k) K(1_{[b_{k}]}, \omega).\nonumber
\end{align}
So, for every $c\in S$ with $c\neq a_i$, 
\begin{align}
\label{eq:KK_ineqaulity}
\\
K(1_{[a_{1},\dots, a_N]}, c) \geq &C_0^{-1}e^{\phi_{N-1}(a_{1},\dots, a_{N}x_{a_N})}K(a_{N},c) \quad(\because \text{ Proposition \ref{prop:green_main_prop}.\ref{prop:green_main_prop_K}})\nonumber \\
 \geq& C_0^{-3}e^{\phi_{N-1}(a_{1},\dots, a_{N}x_{a_N})} \frac{F(a_{N},c)}{F(o,c)}\quad (\because \text{ Proposition \ref{prop:green_main_prop}.\ref{prop:green_main_prop_FG}})\nonumber \\
\geq& C_{0}^{-4}e^{\phi_{N-1}(a_{1},\dots, a_{N}x_{a_N})} \frac{F(a_{N},b_k) F(b_k,c)}{F(o,c)} \quad (\because \text{ Proposition \ref{prop:green_main_prop}.\ref{prop:green_main_prop_F}})\nonumber\\
\geq &C_{0}^{-6} e^{\phi_{N-1}(a_{1},\dots, a_{N}x_{a_N})}F(a_{N},b_k)K(b_{k}, c)  \quad (\because \text{Proposition \ref{prop:green_main_prop}.\ref{prop:green_main_prop_FG} })\nonumber\\
\geq &C_{0}^{-11} C_{3^{}}^{-1}K([a_{1},\dots, a_N],b_{k})\frac{K(b_{k}, c)}{K(1_{[b_{k}]}, \omega)}.  \quad (\because \text{Eq. (\ref{eq:omega_ext_eq})})\nonumber
\end{align}
 Let $\omega' \in \mathcal{M}_m$ and let $c_n\in S$ with $\lim_{n\rightarrow\infty}K(f, c_n)=K(f,\omega')$ for every $f\in C_c(X)$. Such a sequence exists by Corollary \ref{cor:conv_seq_states}.  Since Eq. (\ref{eq:KK_ineqaulity}) holds for all $n$ large enough, 
$$\mu_{\omega'}([a_{1}, \dots, a_N]) \geq   C_{0}^{-11} C_{3^{}}^{-1}K([a_{1},\dots, a_N],b_{k}) \frac{\mu_{\omega'}([b_k])}{\mu_{\omega}([b_k])}. $$
Recall that for every $\mu\in \mathcal{C}$ there exists a finite measure $\nu$ on $\mathcal{M}_m$ s.t. $\mu = \int K(\cdot, \omega')d\nu(\omega')$. Therefore, for every $\mu\in \mathcal{C}$ and every $f\in C_c^+(X)$ we have that
$$\mu([a_{1}, \dots, a_N]) \geq   C_{0}^{-11} C_{3^{}}^{-1}K([a_{1},\dots, a_N],b_{k}) \frac{\mu([b_k])}{\mu_{\omega}([b_k])} $$
Again, since the collection of cylinder sets linearly spans
a dense subset of $C_c(X)$, the above inequality holds for all $f\in C_c^+(X)$.
We take $k\rightarrow\infty$ and obtain that for every $\mu\in \mathcal{C}$ and every $f\in C_c^+(X)$,
\begin{equation}
\label{eq:hyperbolic_main_thm_minimality_ineq}
\mu(f) \geq C_{0}^{-11} C_{3^{}}^{-1} \mu_\omega(f)\limsup_{k\rightarrow\infty} \frac{\mu([b_{k}])}{\mu_\omega([b_{k}])}.
\end{equation}
Let $\mu \in \mathcal{B}_\omega$ and let $\mu' = \mu_\omega-\mu$. By definition of $\mathcal{B}_\omega$, $\mu' \geq 0$ and so $\mu'\in C$. \ 
Since $\inf_{f\in C_c^+(X)}\frac{\mu'(f)}{\mu_\omega(f)}=0$, Eq. (\ref{eq:hyperbolic_main_thm_minimality_ineq}) with $\mu'$ implies that $\lim_{k\rightarrow\infty} \frac{\mu'([b_{k}])}{\mu_\omega([b_{k}])}=0$. In particular, $\lim_{k\rightarrow\infty} \frac{\mu([b_k])}{\mu_\omega([b_k])}=1$. We use this fact and apply Eq. (\ref{eq:hyperbolic_main_thm_minimality_ineq}) with $\mu$ to obtain that $\mu \geq \epsilon_{2} \mu_\omega$ where $\epsilon_2 = C_{0}^{-11} C_{3^{}}^{-1}\in (0,1)$. 

Set $c_n = \epsilon_2 \bigl( 1+ (1-\epsilon_2) + \dots + (1-\epsilon_2)^n\bigr)$. We show by induction that for every $\mu\in \mathcal{B}_\omega$ and every $n\geq 0$,  $\mu\geq c_n\mu_\omega $. Since $c_0 = \epsilon_2$, it is true for $n=0$. Let $n>0$ and suppose that $\mu\geq c_k\mu_\omega $ for every $\mu \in \mathcal{B}_\omega$ and every $k<n$. Then, for every $\mu \in \mathcal{B}_\omega$,  $\frac{1}{1-c_{n-1}}(\mu - c_{n-1}\mu_\omega)\in \mathcal{B}_\omega$  and so $\frac{1}{1-c_{n-1}}(\mu - c_{n-1}\mu_\omega) \geq \epsilon_2 \mu_\omega$. In particular,  $\mu\geq (c_{n-1}+\epsilon_2(1-c_{n-1}))\mu_\omega=c_{n }\mu_\omega$. Letting $n\rightarrow\infty$, we get that $\mu \geq \mu_\omega$. Therefore $\mu = \mu_{\omega}$ for every $\mu \in \mathcal{B}_\omega$, namely $\mathcal{B}_\omega = \{\mu_\omega\}$.

In what follows, let $\omega(\xi)\in \mathcal{M}_m$ be the unique limit point s.t. $K(\cdot, b_n) \rightarrow K(\cdot, \omega(\xi))$, where $b_n\rightarrow\xi$ in the hyperbolic geometry. 
Since $\mathcal{A}_\xi$ contains a single point, $\omega(\xi)$  is well-define. 

By Corollary \ref{cor:conv_seq_states}, the mapping $\omega$ is onto. We show that for two  boundary points $\xi, \eta\in \partial S$ with $\eta \neq \xi$, we have that $K(\cdot, \omega(\xi)) \neq K(\cdot, \omega(\eta))$. 
\begin{lemma}
Let $(S, E)$ be a   $\delta$-hyperbolic graph. Then, for two every  boundary points $\xi, \eta \in \partial (S,E)$ there exists a two-sided infinite geodesic segment $(\dots, a_{-1}, a_0, a_1, \dots)$ s.t. $a_n\xrightarrow[n\rightarrow\infty]{} \xi$ and $a_{-n}\xrightarrow[n\rightarrow\infty]{}\eta$.
\end{lemma}
\begin{proof}
See Lemma 22.15 in \cite{woess_2000}.
\end{proof}
By Proposition \ref{prop:green_main_prop_2}.\ref{prop:green_main_prop_harnack} we have that for every $s\in S$, 
\begin{equation}
\label{eq:ancona_ineq_proof_f}
F(o, s)= C_0^{\pm d_{E}(o, a_0)}F(a_0, s).
\end{equation}
Hence, for every $n,k\geq 0$, with $C_3=C_3(0)$,
\begin{align*}
&\frac{K(a_{k}, a_{n})}{K(a_{k}, a_{-n})}\\ \geq & C_0^{-2}\frac{F(a_k, a_n) F(o, a_{-n})}{F(o, a_n) F(a_k, a_{-n})}\qquad(\because \; \text{Proposition }\ref{prop:green_main_prop}.\ref{prop:green_main_prop_K}) \\
\geq & C_0^{-2-2d_{E}(o,a_0)}\frac{F(a_k, a_n) F(a_{0}, a_{-n})}{F(a_{0}, a_n) F(a_k, a_{-n})}\ (\because \; \text{Eq. }(\ref{eq:ancona_ineq_proof_f}))\\
\geq & C_0^{-2-2d_{E}(o,a_0)}C_3^{-2}\frac{F(a_k, a_n) F(a_{0}, a_{-n})}{F(a_{0}, a_k)F(a_{k}, a_n) F(a_k, a_{0}) F(a_0, a_{-n})} \quad(\because \; \text{Theorem }\ref{thm:main_F_ineq})\\
= & C_0^{-2-2d_{E}(o,a_0)}C_3^{-2}\frac{1}{F(a_0, a_k)F(a_k, a_0)}\\
\geq & C_0^{-4-2d_{E}(o,a_0)}C_3^{-2}\frac{G(a_k, a_k)G(a_0 , a_0)}{G(a_0, a_k)G(a_k, a_0)}\quad(\because \; \text{Proposition }\ref{prop:green_main_prop}.\ref{prop:green_main_prop_FG}) \\
\geq & C_0^{-4-2d_{E}(o,a_0)}C_3^{-2}\frac{1}{G(a_0, a_k)G(a_k, a_0)}.\quad(\because \;G(a,a) \geq 1_a(ax_a) = 1) 
\end{align*}
 Letting $n\rightarrow\infty$, we get that
 $$ \frac{K(1_{[a_k]}, \omega(\xi))}{K(1_{[a_{k}]}, \omega(\eta))} \geq  C_0^{-4-2d_{E}(o,a_o)}C_3^{-2}\frac{1}{G(a_0, a_k)G(a_k, a_0)}. $$
By Proposition \ref{prop:green_main_prop}.\ref{prop:green_main_prop_Gt} and the assumption that $P_G(\phi)<0$, we have that $$\sum_{w\in S}G(a_{0}, w)G(w, a_0) < \infty.$$ Therefore, $G(a_0, a_k)G(a_k, a_0)\rightarrow 0$ as $k\rightarrow \infty$. In particular, there exists $k$ s.t. $K(a_{k}, \omega(\xi)) > K(a_{k}, \omega(\eta))$ and thus $\omega(\xi)\neq \omega(\eta)$. 

\hfill\ensuremath{\square}

%% file: appendix.tex
\section{Additional proofs}

\paragraph{\bf{Proof of Proposition \ref{prop:group_ext_transitive}.}}
By \textbf{(Gen)}, for every $\gamma_1,\gamma_2\in \Gamma_0$ there exists an admissible word $a_1,\dots, a_n\in S_0$ s.t. $\gamma_1 \Gamma = e_{a_1\dots a_n}^{-1}\gamma_2\Gamma$. Thus it suffices to show that for every $a,b\in S_0$ there is an admissible word $a_1,\dots, a_n\in S_0$ s.t. $a, a_1,\dots,a_n, b$ is admissible  and $e_{a,a_1,\dots, a_n}\in \Gamma$. 

Recall that $F_0\subseteq \mathbb{D}$ is a fundamental domain of $\mathbb{D} / \Gamma_0$. Let $F\subseteq \mathbb{D}$ be a fundamental domain of $\mathbb{D} / \Gamma$ with $F_0 \subseteq F$. Let $g^t :T^1 (F)\circlearrowleft$ be the geodesic flow. The geodesic flow on the tangent bundle of a normal cover of a compact hyperbolic surface is known to be topologically transitive, see Theorem 3.8 in  \cite{eberlein_1972}. Therefore there exists $(p_{0},\xi_{0})\in T^1 (F)$ s.t. $\{g^t(p_{0},\xi_{0})\}_{t\in \mathbb{R}}$ is a dense orbit.
We interpret $T^1 (F)$ as $F \times \partial \mathbb{D}$. 
For every $a\in S_0$, let $$B_a = \{(p, \xi) \in T^1 (F) : p\in F_0 \text{ and }\xi \in I_a \}.$$
Since $T^1(F) = \cup_{a}B_a$, for every $a\in S_0$ the set $B_a$ has a non-empty interior and in particular there exists $t_a\in \mathbb{R}$ s.t. $g^{t_a}(p_{0},\xi_{0})\in int(B_a)$. Let $a,b\in S_0$ and let $a_{0,}a_1,\dots, a_n$ be the labels  of the edges of the copies of $F_0$ in $F$ that the geodesic curve between $g^{t_a}(p_{0},\xi_{0})$ and $g^{t_b}(p_{0},\xi_{0})$ intersects. In case the curve passes through a vertex of a copy of $F_0$, we perturb the curve around the vertex, see Figure 5 in \cite{series_1986}. By definition of $B_a$, we have that $a_0=a$. Then, we have that $a,a_1,\dots, a_n, b$ is an admissible word and since the curve starts in $F_0$ and finishes in $F_0$,  $e_{a, a_1,\dots, a_n}\in \Gamma$. 
\hfill\ensuremath{\square}

\paragraph{\bf{Proof of Proposition \ref{prop:pressure_delta}.}}
The following arguments are taken from \cite{stadlbauer_2013}. For every $a\in S_{0}$, we fix  $\sigma_{a}\in T_{\Sigma}[a] \subseteq \Sigma$ arbitrarily. Series showed that there exists $C > 1$ s.t. for every admissible $a_1,\dots,a_n\in S_{0}$
\begin{equation}
\label{eq:potential_propto_d}
\exp \bigl(\phi^{\Sigma,\delta }_n(a_1,\dots, a_n \sigma_{a_n})\bigr)  = C^{\pm 1} \exp \bigl (-\delta d_{\mathbb{D}}(0, \gamma ^{-1}0)\bigr) 
\end{equation}  
where $\gamma = e_{a_n}^{-1}, \dots, e_{a_1}^{-1}$.
See Section 5 in \cite{series_1981}. 
Then, by \textbf{(Bnd)},  
\begin{equation}
\label{eq:poincare_green}
\sum_{\gamma \in \Gamma}e^{-\delta d_{\mathbb{D}}(0, \gamma 0)} \asymp \sum_{n=0}^\infty \sum_{\substack{a_1,\dots, a_n \in S_{0}\\ e_{a_1}\cdots e_{a_n} \Gamma = \Gamma}}e^{\phi_{n}^{\Sigma, \delta}(a_1, \dots, a_{n}\sigma_{a_n})} \asymp \sum_{a,b\in S_{0}}\sum_{n=0}^\infty (L_{
\phi^{X,\delta}}^n (1_{[a]\times \{\Gamma\}}))(x_b) 
\end{equation}
with $x_b \in T[(b,\Gamma)]\subseteq X$.
    In particular, the Poincar\'e series and the Green's function converge or diverge together. 

Next, assume $\delta > \delta(\Gamma)$. We show that there exists some $t \in (0,1)$ s.t.
$$ \sum_{n=0}^\infty \sum_{\substack{a_1,\dots, a_n \in S_{0}\\ e_{a_1}\cdots e_{a_n} \Gamma = \Gamma}}t^{-n}e^{\phi_{n}^{\Sigma, \delta}(a_1, \dots, a_{n}\sigma_{a_n})} < \infty. $$
  Let $\epsilon>0$ so that $(1-\epsilon)\delta > \delta(\Gamma)$. By (\ref{eq:potential_propto_d}) there exist $n_0\geq 0$ and $\alpha >0$  s.t. for all $n \geq n_0$,
$$\min_{\sigma\in \Sigma} \{e^{-\epsilon\phi_{n}^{\Sigma,\delta}(\sigma)}\}>e^{-\alpha}.$$ 
  Let $$\beta = \max_{n \leq n_0, \sigma \in \Sigma}|\phi_n^{\Sigma, \delta}(\sigma)| $$ and let $n_1>n_0$  large enough s.t. $\lfloor \frac{n}{n_0}\rfloor \cdot \frac{1}{n}>\frac{1}{2}$ for all $n\geq n_1$. Then, 
\begin{align*}\infty >& \sum_{n=0}^\infty \sum_{\substack{a_1,\dots, a_n \in S_{0}\\ e_{a_1}\cdots e_{a_n}  \Gamma = \Gamma}}e^{\phi_{n}^{\Sigma, (1-\epsilon)\delta}(a_1, \dots, a_{n}\sigma_{a_n})}\\
\geq &
\sum_{n=n_1}^{\infty} \sum_{\substack{a_1,\dots, a_n \in S_{0}\\ e_{a_1}\cdots e_{a_n} \Gamma = \Gamma}}e^{\phi_{n}^{\Sigma,\delta}(a_1, \dots, a_{n}\sigma_{a_n})}  e^{-\epsilon  \phi_{n}^{\Sigma, \delta}(\sigma) }\\
 \geq &
\sum_{n=n_1}^{\infty} \sum_{\substack{a_1,\dots, a_n \in S_{0}\\ e_{a_1}\cdots e_{a_n} \Gamma = \Gamma}}e^{\phi_{n}^{\Sigma,\delta}(a_1, \dots, a_{n}\sigma_{a_n})}  e^{-\epsilon  n_{0}\beta - \lfloor \frac{n}{n_0}\rfloor \epsilon\alpha   }\\
\geq &\ \sum_{n=n_1}^{\infty} \sum_{\substack{a_1,\dots, a_n \in S_{0}\\ e_{a_1}\cdots e_{a_n} \Gamma = \Gamma}}e^{\phi_{n}^{\Sigma,\delta}(a_1, \dots, a_{n}\sigma_{a_n})}  e^{-\epsilon  n_{0}\beta}\bigl(e^{ -\frac{1}{2}\epsilon  \alpha   }\bigr)^n. 
\end{align*}
\hfill\ensuremath{\square}

\paragraph{\bf{Proof of Proposition \ref{prop:non_atomic}.}}
The proof of the proposition relies on the following lemma:
\begin{lemma}
Let $\sigma\in \Sigma$. Then, there exists a sequence of sequences $\{\sigma^{n} \}\subseteq \Sigma$ and $\{\gamma_n\} \subseteq\ \Gamma$ s.t.
\begin{enumerate}
\item 
$\{\sigma^n\}$ and $\{\pi_{\Sigma}(\sigma^n)\}$ are both infinite sets.
\item
$\pi_{\Sigma}(\sigma^n) = \gamma_n \pi_{\Sigma}(\sigma)$.
\item
There exists $B>0$ s.t. for all $n\geq 0$,\ $|\gamma_n'(\pi_\Sigma(\sigma^n))|>B$. 
\item For every $n\geq 0$ there exist $m_n, k_n\geq0$ s.t. for every $\gamma\in \Gamma_0$,  $$T^{m_n}(\sigma^n,\gamma \Gamma) = T^{k_n}(\sigma,\gamma \Gamma)$$ and $$\exp \biggl(\phi^{\Sigma,\delta}_{m_n}(\sigma^n)-\phi^{\Sigma,\delta}_{k_n}(\sigma) \biggr)= |\gamma_n'(\pi_\Sigma(\sigma_n))|^{\delta}.$$ 
\end{enumerate}
\end{lemma}
\begin{proof}
For every $ a,b\in S_{0}$ let $w_{a,b}$ be an admissible word that  includes $b$ and such that $aw_{a,b} a$ is an admissible word. By Proposition \ref{prop:group_ext_transitive}, we can choose $w_{a,b}$ so that $e_{aw_{a,b}}\in \Gamma$.
We also choose $w_{a,b}$ so that for every $a_1,a_2\in S_{0}$ with $a_1\neq a_2$, we have that $w_{a_1,b}\neq w_{a_2,b}$. 

Let $\xi = \pi_{\Sigma}(\sigma)$ and let$$\sigma^{n,b} = (\sigma_0, \dots, \sigma_{n-1} ,\sigma_{n} w_{\sigma_{n},b }\sigma_n,\sigma_{n+1},\dots) .$$  
For every $n$, let $b_n \in S_{0}$ s.t. $\pi_{\Sigma}(\sigma^{n,b_n})\neq \xi$, let $\sigma^n = \sigma^{n,b_n}$ and let $\xi_n = \pi_\Sigma(\sigma^n)$. Such $b_n$ exists since $|S_{0}|\geq 4$ and $\pi_{\Sigma}$ is at most two-to-one (see \cite{series_1986}). Let $\alpha_n\in \Gamma_0$ s.t. $f_{\Gamma_0}^n \xi = \alpha_n^{-1}\xi$. By \textbf{(Res)}, for all $n$, $f_{\Gamma_0}^n \xi_n = \alpha_n^{-1}\xi_n$. Let $l_n = |\sigma_{n} w_{\sigma_{n},b_n}|$ and let $\beta_n\in \Gamma_0$ s.t. $f_{\Gamma_0}^{l_n}(f_{\Gamma_0}^n \xi_n) = \beta_n (f_{\Gamma_0}^n \xi_n)$. 
Since $f_{\Gamma_0}^{n}\xi = f_{\Gamma_0}^{n+l_n}\xi_n$, we have that 
$$\xi_n = \alpha_n \beta_n^{-1} \alpha_n^{-1} \xi.  $$
Let $\gamma_n = \alpha_n \beta_n^{-1} \alpha_n^{-1} $. Observe that  $\beta_n = e^{-1}_{\sigma_n w_{\sigma_n,b_{n}}}\in \Gamma$ and thus $\gamma_n \in \Gamma$. In particular, for all $\gamma\in \Gamma$, $$T^{n}(\sigma, \gamma\Gamma) = T^{n + l_n}(\sigma^n, \gamma \Gamma).$$
 We deduce,
\begin{align*}
 |\gamma_n'(\xi)| = &  | \alpha_n'(\beta_n^{-1} \alpha_n^{-1}\xi)| \cdot| (\beta_n^{-1})'(\alpha_n^{-1}\xi)| \cdot |(\alpha_n^{-1})'(\xi)|\\
=&| (\beta_n^{-1})'(\alpha_n^{-1}\xi)| \cdot\frac{   |(\alpha_n^{-1})'(\xi)| }{|(\alpha_n^{-1})'(\alpha_n\beta_n^{-1} \alpha_n^{-1}\xi)|}  \\
=& | (\beta_n^{-1})'(f_{\Gamma_0}^n \xi)| \cdot\frac{   |(f_{\Gamma_0}^{n})'(\xi)| }{|(f_{\Gamma_0}^{n})'(\xi_{n})|}.
\end{align*}
By \textbf{(Dist)}, there exists $B>1$ s.t. 
$$\frac{   |(f_{\Gamma_0}^{n})'(\xi)| }{|(f_{\Gamma_0}^{n})'(\xi_{n})|} \geq B^{-1}.$$
 Then, with
$$D = \min_{a,b\in S_{0}} \min_{\xi'\in [a]} |(e_{a w_{a,b}}^{-1})'(\xi')| > 0$$ 
we have that
$$|\gamma_n'(\xi)|  \geq \frac{D}{B}. $$

We show that $|\{\xi_n\}|=|\{\sigma^n\}|=\infty$. The mapping $\pi$ is continuous and thus $\xi_n \rightarrow \xi$. Since $\xi_n \neq \xi$ for all $n$, there exists a sub-sequence $\xi_{n_k}$ with $|\{\xi_{n_k}\}|=\infty$. 

Lastly, since
$$\exp \bigl(\phi_{n}^{\Sigma, \delta}(\sigma)\bigr) = |(\alpha^{-1}_{n})'(\xi)|^{-\delta}  $$     
and
$$\exp \bigl(\phi_{n +l_{n}}^{\Sigma, \delta}(\sigma^n)\bigr) = |(\beta_{n}\alpha_{n}^{-1})'(\xi_{n})|^{-\delta}= |(\alpha_n \beta_n^{-1})'(\alpha_n ^{-1}\xi)|^\delta $$
the lemma follows with $k_n = n$ and $m_n = n + l_n$. 
\end{proof}
\begin{enumerate}[label={(\arabic*)}]
\item
Assume by contradiction that $\xi\in \mathbb{D}$ is an atom. By the auxiliary lemma, with $\sigma \in \pi^{-1}_\Sigma(\xi)$, there exists a sequence $\gamma_n\in \Gamma$ s.t. $\{\gamma_n \xi\}$ are all distinct and  $|\gamma_n'(\xi)|^\delta$ is bounded from below. Then, 
$$\mu\bigl(\{\gamma_n \xi\}_{n\in \mathbb{N}}\bigr)  =\sum_{n}\mu\bigl(\{\gamma_n \xi\}\bigr) =\sum_n |\gamma_n'(\xi)|^{\delta } \mu(\{\xi\}) = \infty$$
which contradicts the finiteness of $\mu$.

\item

Let $\mu$ be a Radon measure with $L_{\phi^{X,\delta}}^*\mu = \mu$. Observe that 
if $x,y\in X$ with $T^nx = T^my$, then 
$$e^{-\phi_n^{X,\delta}(x)}\mu(\{x\}) = e^{ - \phi^{X,\delta}_m(y)}\mu(\{y\}). $$
This is because 
$$\mu(\{x\}) = \mu(L_{\phi^{X,\delta}}^n \delta_x) = e^{\phi_n^{X,\delta}(x)}\mu(\{T^nx\}) $$
and
$$\mu(\{y\}) = \mu(L_{\phi^{X,\delta}}^m \delta_y) = e^{\phi_m^{X,\delta}(y)}\mu(\{T^{m}y\}). $$
Assume by contradiction that $(\sigma, \gamma\Gamma)\in X$ is an atom. By the auxiliary lemma, there exists a sequences $\sigma^n\in \Sigma$ and $m_n, k_n\geq 0$ s.t. $\{\sigma^n\}$ are all distinct,  $T^{m_n}(\sigma^n, \gamma\Gamma) = T^{k_n}(\sigma, \gamma\Gamma)$  for every $\gamma\in \Gamma_0$ and  $\phi^{\Sigma,\delta}_{m_n}(\sigma^n)-\phi^{\Sigma,\delta}_{k_n}(\sigma)$ is bounded from below. Then, for every $\gamma\in \Gamma_0$,
$$\mu\bigl(\{(\sigma^n, \gamma\Gamma)\}_{n\in \mathbb{N}}\bigr)  =\sum_n \exp \bigl( \phi^{\Sigma,\delta}_{m_n}(\sigma^n)-\phi^{\Sigma,\delta}_{k_n}(\sigma)   \bigr) \mu(\{x\}) = \infty$$
which contradicts the fact that  $\mu$ is a Radon measure.
\end{enumerate}
\hfill\ensuremath{\square}